\documentclass[11pt, a4paper]{amsart}

\usepackage[T1]{fontenc}
\usepackage[utf8]{inputenc}
\usepackage[english]{babel}
\usepackage[leqno]{amsmath}
\usepackage{amssymb}
\usepackage{tabularx}
\usepackage{a4wide}
\usepackage{url}
\usepackage{tikz-cd}
\usepackage{csquotes}
\usepackage{enumitem}
\usepackage{mathtools}
\usepackage{manfnt}
\usepackage{mathrsfs}
\usepackage{aurical}
\usepackage[leqno]{amsmath}

\usepackage[
  maxnames=4,
  maxalphanames=4,
  style=alphabetic,
  url = false,
  doi = false,
  isbn = false,
  ]{biblatex}
\renewbibmacro{in:}{}

\usepackage[disable]{todonotes} 

\definecolor{cadmiumgreen}{rgb}{0.0, 0.42, 0.24}
\definecolor{darkred}{rgb}{.85,0,0}
\usepackage[
  colorlinks,
  linkcolor=darkred,
  urlcolor=darkred,
  citecolor=cadmiumgreen,
  pdfstartview ={FitV},
]{hyperref}
\hypersetup{
  pdfauthor={Omid Amini, Matthieu Piquerez},
  pdftitle={Homological smoothness and Deligne resolution for tropical fans},
  bookmarksdepth = 2}

\usepackage{scalerel}

\usetikzlibrary{calc}
\tikzset{vcenter/.style={baseline={([yshift=-.8ex]current bounding box.center)}}}
\tikzset{dot/.style={insert path={node {\tikz[baseline=.6pt]\filldraw[black] (0,0) circle (1.2pt);}}}}

\setlist[itemize,1]{itemsep=\smallskipamount}
\setlist[enumerate,1]{itemsep=\smallskipamount, label=\textnormal{(\arabic*)}}

\newtheorem{thm}{Theorem}[section]

\newtheorem{prop}[thm]{Proposition}

\theoremstyle{definition}

\newtheorem{remm}[thm]{Remark}
\newenvironment{remark}
  {\pushQED{\qed}\remm}
  {\popQED\endremm}
\newtheorem{exx}[thm]{Example}
\newenvironment{example}
  {\pushQED{\qed}\exx}
  {\popQED\endexx}

\numberwithin{equation}{section}

\newcommand{\ie}{i.e.}

\newcommand{\resp}{resp.\ }

\renewcommand{\~}{\widetilde}
\renewcommand{\hat}{\widehat}
\newcommand{\Q}{\mathbb{Q}}
\newcommand{\Z}{\mathbb{Z}}

\newcommand{\R}{\mathbb{R}}

\newcommand{\vsim}{\rotatebox{90}{$\sim$}}
\newcommand{\longhookrightarrow}{\lhook\joinrel\longrightarrow}

\let\oldforall\forall
\renewcommand{\forall}{\oldforall\:}
\let\oldbigwedge\bigwedge
\renewcommand{\bigwedge}{{\textstyle\oldbigwedge\!}}

\renewcommand{\geq}{\geqslant}
\renewcommand{\leq}{\leqslant}
\renewcommand{\setminus}{\smallsetminus}

\makeatletter
\let\oldsum\sum
\renewcommand{\sum}{\@ifnextchar_\@mysum\oldsum}
\def\@mysum_#1{\oldsum_{\substack{#1}}}
\let\oldbigoplus\bigoplus
\renewcommand{\bigoplus}{\@ifnextchar_\@mybigoplus\oldbigoplus}
\def\@mybigoplus_#1{\oldbigoplus_{\substack{#1}}}
\let\oldprod\prod
\renewcommand{\prod}{\@ifnextchar_\@myprod\oldprod}
\def\@myprod_#1{\oldprod_{\substack{#1}}}
\makeatother

\makeatletter
\let\oldnu\nu
\newlength{\heightnu}
\newlength{\depthnu}
\def\nu#1_#2{{\settoheight{\heightnu}{\hbox{$#2$}}\settodepth{\depthnu}{\hbox{$#2$}}\oldnu\rule[\depthnu-3pt]{0pt}{1pt}#1_{\!#2}}}
\makeatother

\newcommand{\ssnu}{\ssub{\oldnu}!}

\makeatletter
\let\@oldinfty\infty
\newcommand{\@sminfty}{{\scaleto{\@oldinfty}{2.8pt}}} 
\renewcommand{\infty}{{\mathchoice%
  {\displaystyle{\@oldinfty}}%
  {\textstyle{\@oldinfty}}%
  {\scriptstyle{\@sminfty}}%
  {\scriptscriptstyle{\@sminfty}}}
}
\makeatother

\newcommand{\rquot}[2]{#1\big/#2}
\newcommand{\rest}[1]{\raisebox{-1pt}{$\vert$}_{#1}}

\newcommand{\id}{\mathrm{id}} 
\newcommand{\dual}{\star}
\newcommand{\st}{\bigm|} 
\newcommand{\bul}{\bullet}

\DeclareMathOperator{\rk}{rk} 
\let\hom\relax
\DeclareMathOperator{\hom}{Hom} 
\let\Im\relax
\DeclareMathOperator{\Im}{Im} 
\DeclareMathOperator{\ord}{ord} 
\DeclareMathOperator{\coker}{coker}

\renewcommand{\i}{\mathrm i} 
\renewcommand{\d}{\mathrm d} 
\newcommand{\hyp}{\mathbb H} 
\newcommand{\BM}{{^{\scaleto{\mathrm{BM}}{3.6pt}}}}
\DeclareMathOperator{\Cone}{Cone} 
\DeclareMathOperator{\sign}{sign} 

\newcommand{\T}{\mathbb T} 
\renewcommand{\P}{\mathbb P} 
\newcommand{\TP}{{\T\P}} 

\newcommand{\SF}{\ss[-1pt][0pt]{\mathbf F}} 
\newcommand{\RMod}{{\R_+}} 
\newcommand{\conezero}{{\underline0}}
\newcommand{\e}{\ss{\mathfrak e}} 
\newcommand{\nvect}{\ss{\mathfrak n}} 
\newcommand{\chart}{\~} 
\newcommand{\cube}{\ssbis[0pt][.5pt]{\scaleto{\square}{6pt}}!} 
\DeclareMathOperator{\sed}{sed} 

\newcommand{\shiftcomp}[2][0]{{}\mkern#1mu\overline{\mkern-#1mu#2}}
\newcommand{\comp}[1]{\if#1X \shiftcomp[3]{#1}\else\if#1Z \shiftcomp[3]{#1} \else \shiftcomp{#1}\fi\fi} 

\newcommand{\suppaux}[2]{\scalebox{1}[1.4]{$#1\lvert$}#2\scalebox{1}[1.4]{$#1\rvert$}}
  \newcommand{\supp}[1]{\mathpalette\suppaux{#1}}
\newcommand{\dimsaux}[2]{\raisebox{.2ex}{\scalebox{1}[.8]{$#1\lvert$}}#2\raisebox{.2ex}{\scalebox{1}[.8]{$#1\rvert$}}}
  \newcommand{\dims}[1]{\mathpalette\dimsaux{#1}}
\newcommand{\subface}{\prec}
\newcommand{\ssubface}{\mathbin{\mathchoice
  {\subface\!\!\!\cdot}%
  {\subface\!\!\!\cdot}%
  {\subface\!\cdot}%
  {\subface\!\cdot}%
}} 
\newcommand{\supface}{\succ}
\newcommand{\ssupface}{\mathbin{\mathchoice
  {\cdot\!\!\!\supface}%
  {\cdot\!\!\!\supface}%
  {\cdot\!\supface}%
  {\cdot\!\supface}%
}}

\newcommand{\LL}{\mathscr L} 
\newcommand{\eLL}{\hat\LL} 

\renewcommand{\div}{\mathrm{div}} 
\newcommand{\Div}{\mathrm{Div}} 
\newcommand{\Prin}{\mathrm{Prin}} 

\newcommand{\Csh}{\mathscr C} 
\newcommand{\Bsh}{\mathscr B} 
\newcommand{\Ssh}{\mathscr S} 
\newcommand{\Sh}{\Ssh} 

\newcommand{\tropmod}[2]{{\mathcal{T\!M}}_{\!#1}(#2)} 
\newcommand{\basetm}[1]{#1_o} 
\newcommand{\symbuptm}{\tikz[scale=.2, baseline=-.1]{\draw(.1,-.05)to[out=50,in=-130](.9,.05) (.5,0)--++(0,.5);}}
\newcommand{\uptm}[1]{{#1}_{\!\symbuptm}} 

\newcommand{\etm}{\e_{\!\symbuptm}} 
\newcommand{\prtm}{\mathfrak p} 

\DeclareMathOperator{\gys}{Gys} 

\newcommand{\Ma}{{\scalebox{1.12}{$\mathfrak m$}}} 
\newcommand{\distel}{\ast} 

\RenewDocumentCommand{\ss}{O{0pt} O{0pt} O{.9} m e{_^}}{
  #4%
  \IfValueT{#5}{
    \sb{\hspace{#1}\scaleobj{#3}{#5}}
  }
  \IfValueT{#6}{
    \sp{\hspace{#2}\scaleobj{#3}{#6}}
  }
}

\newcommand{\E}{\textnormal{\textsf{E}}}

\newcommand{\vv}{\mathbf v} 
\newcommand{\cycl}{\mathrm{cl}} 
\newcommand{\ones}{(1, \dots, 1)}

\usepackage{bbding}

\newcommand{\longsimto}{\xrightarrow{\ \raisebox{-3pt}[0pt][0pt]{\small$\hspace{-1pt}\sim$\ }}}

\newcommand{\mT}{\mathcal T\!} 

\NewDocumentCommand{\ssub}{O{0pt} O{.9} m t! e{_^}}{
  #3%
  \IfValueT{#5}{
    \IfBooleanTF{#4}{\sb{\hspace{#1}\scaleobj{#2}{#5}}}{\sb{#5}}
  }
  \IfValueT{#6}{
    \IfBooleanTF{#4}{\sp{\hspace{#1}\scaleobj{#2}{#6}}}{\sp{#6}}
  }
}

\NewDocumentCommand{\ssbis}{O{0pt} O{0pt} O{.8} m t! e{_^}}{
  #4%
  \IfValueT{#6}{
    \IfBooleanTF{#5}{\sb{\hspace{#1}\scaleobj{#3}{#6}}}{\sb{#6}}
  }
  \IfValueT{#7}{
  \IfBooleanTF{#5}{\sp{\hspace{#2}\scaleobj{#3}{#7}}}{\sp{#7}}
}
}
\NewDocumentCommand{\ssubis}{O{0pt} O{0pt} O{.9} m t! e{_^}}{
  #4%
  \IfValueT{#6}{
    \IfBooleanTF{#5}{\sb{\hspace{#1}\scaleobj{#3}{#6}}}{\sb{#6}}
  }
  \IfValueT{#7}{
  \IfBooleanTF{#5}{\sp{\hspace{#2}\scaleobj{#3}{#7}}}{\sp{#7}}
}
}
\ExplSyntaxOn
\NewDocumentCommand{\tossub}{o o m}{
  \expandafter\let\csname old\cs_to_str:N #3\endcsname#3
  \renewcommand#3%
  {\ssub[#1][#2]{\csname old\cs_to_str:N #3\endcsname}}
}
\ExplSyntaxOff

\newcommand{\ssinfty}{\ssub{\infty}!}

\tossub\omega
\tossub\varpi
\tossub\deg
\tossub\LL
\tossub\eLL
\tossub\e
\tossub\distel
\tossub\rk
\tossub\P
\tossub[-2pt][1]\Gamma
\tossub\Div
\tossub\Prin
\tossub\cycl
\tossub\TP
\tossub\i
\tossub\gys
\newcommand{\sssigma}{\ssubis[-2pt]{\sigma}!}
\newcommand{\ssM}{\ssub{M}!}
\newcommand{\ssN}{\ssub{N}!}

\newcommand{\ssSigma}{\ssubis{\Sigma}!}
\newcommand{\compSigma}{\ssubis{\comp\Sigma}!}
\newcommand{\compDelta}{\ssubis{\comp\Delta}!}

\newcommand{\sspi}{\ssub{\pi}!}

\newcommand{\ssf}{\ssub{f}!}
\newcommand{\sseta}{\ss{\eta}}

\newcommand{\ssTP}{\ss\TP}
\newcommand{\ssH}{\ss{H}}
\newcommand{\sse}{\ss\e}
\newcommand{\ssvarpi}{\ss{\varpi}}
\newcommand{\ssi}{\ss\i}
\newcommand{\contract}{\kappa}
\newcommand{\sstau}{\ss{\tau}}
\newcommand{\ssomega}{\ss\omega}

\newcommand{\gst}[2]{\ssub{\langle#1\rangle}!_{#2}} 
\newcommand{\Bsho}{\ssub{\mathcal B}!} 

\def\roundface#1#2[#3]{\fill[#3] (I) -- (I#1) .. controls ($.3*(I#1)+.7*(I#1#1#2)$) and ($.3*(I#1#2)+.7*(I#1#1#2)$) .. (I#1#2) .. controls ($.3*(I#1#2)+.7*(I#2#1#2)$) and ($.3*(I#2)+.7*(I#2#1#2)$) .. (I#2) -- cycle}

\makeatletter
\@namedef{subjclassname@2020}{\textup{2020} Mathematics Subject Classification}
\makeatother

\newcommand{\msc}[1]{\href{http://www.ams.org/msc/msc2020.html?t=&s=#1}{#1}}


\begin{document}

\allowdisplaybreaks

\title{Homological smoothness and Deligne resolution for tropical fans}

\author{Omid Amini}
\address{CNRS - CMLS, \'Ecole polytechnique, Institut polytechnique de Paris.}
\email{\href{omid.amini@polytechnique.edu}{omid.amini@polytechnique.edu}}

\author{Matthieu Piquerez}
\address{LS2N, Inria, Nantes Université}
\email{\href{matthieu.piquerez@univ-nantes.fr}{matthieu.piquerez@univ-nantes.fr}}

\date{\today}

\begin{abstract} We say that a tropical fan is \emph{homologically smooth} if each of its open subsets verify tropical Poincaré duality. A \emph{tropical homology manifold} is a tropical variety that is locally modelled by open subsets of homologically smooth tropical fans.

We show that homological smoothness is a $\mT$-stable property in the category of tropical fans. This implies in particular that quasilinear fans are homologically smooth, and tropical varieties locally modelled by them are tropical homology manifolds. Previously, this was known only for locally matroidal tropical varieties.

In order to show the above results, we prove a tropical analogue of the Deligne weight spectral sequence for homologically smooth tropical fans. This allows to describe the cohomology of tropical modifications, and will be of importance in our companion work which develops a Hodge theory in the tropical setting.
\end{abstract}

\maketitle

\setcounter{tocdepth}{1}

\tableofcontents

\section{Introduction} \label{sec:intro}

The work of this paper is motivated by the development of a Hodge theory in the combinatorial and tropical setting. We study a tropical notion of homology manifolds in the spirit of the classical notion of homology manifolds in topology~\cites{Mio00, Wei02}. As in the classical setting, this will be a local notion. Since tropical fans and their support fanfolds are building blocks of more general tropical varieties, this is reduced to considering tropical fans and tropical fanfolds.

\subsection{Tropical fans}

Let $N \simeq \Z^n$ be a lattice of finite rank and denote by $\ssN_\R$ the vector space generated by $N$. Consider a rational fan $\Sigma$ in $\ssN_\R$. We denote the support of $\Sigma$ by $\supp\Sigma$ and call it a \emph{(rational) fanfold} in $\ssN_\R$. We say that $\Sigma$ is \emph{pure dimensional} if all its maximal cones have the same dimension. For each cone $\sigma \in \Sigma$, denote by $\ssN_\sigma$ the sublattice of $N$ generated by the integer points of $\sigma$.

A \emph{tropical orientation} of a rational fan $\Sigma$ of pure dimension $d$ is an integer valued map
\[ \omega \colon \Sigma_d \to \Z\setminus\{0\} \]
that verifies the \emph{balancing condition}: for each cone $\tau$ in $\Sigma$ of codimension one, we have
\[ \sum_{\sigma \supset \tau} \omega(\sigma)\sse_{\sigma}^{\tau} =0 \qquad \textrm{ in } \rquot{N}{\ssN_\tau} \]
where the sum is over facets $\sigma$ of $\Sigma$ which contain $\tau$, and $\sse_{\sigma}^{\tau}$ is the generator of the quotient $\rquot{(\sigma \cap N)}{(\tau \cap N)} \simeq \Z_{\geq 0}$.

A \emph{tropical fan} in $\ssN_\R$ is a pair $(\Sigma, \omega!_{\Sigma})$ consisting of a pure dimensional rational fan $\Sigma$ endowed with a tropical orientation $\omega!_{\Sigma}$. We refer to the support of $(\Sigma, \omega!_{\Sigma})$ as a \emph{tropical fanfold}.

\subsection{Homological smoothness and tropical homology manifolds}

For a rational fan $\Sigma$, and open subset $U$ of $\supp\Sigma$, let $H^{p,q}(U)$, $H^{p,q}_c(U)$, $H_{p,q}(U),$ and $H^\BM_{p,q}(U)$ be the tropical cohomology, tropical cohomology with compact support, tropical homology, and Borel-Moore homology groups of $U$ with rational coefficients, respectively, introduced in~\cite{IKMZ}. We recall the relevant definitions in Section~\ref{sec:prel}.

Let $\Sigma$ be a tropical fan of dimension $d$. The Borel-Moore homology group $H_{d,d}^\BM(\Sigma)$ contains a canonical element $\ssnu_\Sigma$, the analogue of the fundamental cycle in the tropical setting. Using the cap product $\frown$, we get a natural map
\[ \cdot\frown\ssnu_{\Sigma}\colon H^{p,q}(\Sigma) \longrightarrow H^\BM_{d-p,d-q}(\Sigma). \]

We say that a tropical fan $\Sigma$ verifies \emph{tropical Poincaré duality} if the above map is an isomorphism for each bidegree $(p,q)$.

More generally, for each open subset $U$ of $\supp\Sigma$,
we get by restriction a canonical element $\ssnu_{U}$ in the Borel-Moore homology $H^\BM_{d,d}(U)$. This leads to the natural map
\[
  \cdot\frown\ssnu_{U}\colon H^{p,q}(U) \longrightarrow H^\BM_{d-p,d-q}(U).
\]
We say that $\Sigma$ is \emph{homologically smooth}, or that $\Sigma$ is a \emph{tropical homology manifold}, if these maps are all isomorphisms, for open subsets $U$ of $\supp{\Sigma}$. Equivalently, the canonical element gives a map $H^{d,d}_c(U) \to \Q$, and tropical Poincaré duality is the statement that the natural pairing
\[ H^{p,q}(U) \times H^{d-p,d-q}_c(U) \to H^{d,d}_c(U) \to \Q \]
is perfect for each open subset $U \subseteq \supp\Sigma$ and all pairs of integers $p,q$.

By definition, homological smoothness is a local condition, and it is a property of the support, that is, it is the property of the underlying tropical fanfold. Moreover, it follows from the Künneth formula that homological smoothness is closed under products.

We will show the following result in Section~\ref{sec:prel}.

\begin{thm}\label{thm:tropical_homology_manifold_local_fans} A tropical fan $\Sigma$ is homologically smooth if and only if\/ $\Sigma$ and all its star fans $\ssSigma^\sigma$, $\sigma \in \Sigma$, verify tropical Poincaré duality.
\end{thm}

In this paper, the star fan $\ssSigma^\sigma$ refers to the fan in $\rquot {\ssN_\R}{\ssN_{\sigma,\R}}$ induced by the cones $\eta$ in $\Sigma$ which contain $\sigma$ as a face.

\smallskip
We say that a tropical variety $X$ is a \emph{tropical homology manifold} if it is locally modelled by tropical fanfolds which are homologically smooth, in the sense that each point $x$ in $X$ admits an open neighborhood isomorphic to an open subset of $\T^k \times \supp{\Sigma}$ for $\T = \R \cup\{+\infty\}$, $k$ a non-negative integer, and $\Sigma$ a homologically smooth tropical fan. It follows from Theorem~\ref{thm:tropical_homology_manifold_local_fans} that this is equivalent to requiring that each local fanfold of $X$ verify tropical Poincaré duality.

 \smallskip
Local duality implies the following global duality theorem.

\begin{thm}[Poincaré duality for tropical homology manifolds] \label{thm:PD-intro}
  Let $X$ be a tropical homology manifold of dimension $d$. The cap product with the fundamental class $\ssnu_{X}\in H_{d,d}^\BM(X)$ induces an isomorphism
   \[ \cdot \frown \ssnu_{X}\colon H^{p,q}(X) \to H_{d-p,d-q}^\BM(X). \]
  In particular, we have $H^{d,d}_c(X) \simeq \Q$ and the natural pairing
   \[ H^{p,q}(X) \times H^{d-p,d-q}_c(X) \to H^{d,d}_c(X) \simeq \Q \]
  is perfect.
\end{thm}

To a matroid $\Ma$ with ground set $E$, one associates a unimodular fan $\ssSigma_\Ma$ called the \emph{Bergman fan of $\Ma$} in the real vector space $\rquot{\R^E\!}{\R\ones}$. Its support is called a \emph{Bergman fanfold}. It follows from the calculation of homology groups of Bergman fanfolds in~\cites{Sha13a, JSS19}, and the observation that star fans of Bergman fans are themselves Bergman, that a Bergman fanfold is homologically smooth. The proof of Theorem~\ref{thm:PD-intro} is similar to the treatment of Poincaré duality for tropical varieties modeled locally by Bergman fanfolds~\cites{JSS19, GS19}. Our Theorem~\ref{thm:smooth_shellable_intro} below produces a large class of tropical fans that are homologically smooth.

\subsection{$\mT$-stability of homological smoothness}

We prove the following  stability theorem.

\begin{thm} \label{thm:smooth_shellable_intro}
  Homological smoothness is $\mT$-stable. In particular, any quasilinear tropical fan is homologically smooth.
\end{thm}

The notion of $\mT$-stability is introduced in~\cite{AP-hodge-fan}, we will review the definition in Section~\ref{subsec:T_stability}. We recall here that a tropical fan is quasilinear if it is either the real line with an arbitrary tropical orientation, or, it is a product of two quasilinear fans, or, it is the stellar subdivision or stellar assembly of a quasilinear fan, or, it is obtained as a result of tropical modification of a quasilinear fan along a tropical divisor which is itself quasilinear.

\smallskip
The proof of the above result boils down to the following theorem. We will recall the definition of tropical modification in Section~\ref{subsec:tropical_modification}.

\begin{thm} \label{thm:tropical_modification_intro}
  Homological smoothness is preserved under tropical modifications along divisors which are homologically smooth.
\end{thm}

The proof of this theorem is based on explicit description given in Section~\ref{sec:homology_tropical_modification} of the cohomology groups of tropical modifications, that we hope to be of independent interest. This uses the tropical Deligne resolution Theorem~\ref{thm:deligne-intro}, see Section~\ref{sec:Deligne-intro}.

\smallskip
The motivation behind the above results comes from the development of a Hodge theory for tropical varieties. Consider an effective tropical fan $(\Sigma, \omega!_\Sigma)$, that is, a tropical fan whose orientation $\omega!_\Sigma$ takes only positive values. We say that $\Sigma$ is Kähler if it is homologically smooth and it is Chow-Kähler in the sense of~\cite{AP-hodge-fan}, that is, it is quasi-projective, and the Chow ring of each star fan of $\Sigma$ verifies Hard Lefschetz property and Hodge-Riemann bilinear relations. Combining Theorem~\ref{thm:smooth_shellable_intro} with $\mT$-stability of being Chow-Kähler, established in~\cite[Theorem 1.9]{AP-hodge-fan}, we obtain the following result.

\begin{thm}
  Being Kähler is $\mT$-stable in the class of  quasi-projective tropical fans which are effective and simplicial. In particular, a  simplicial quasi-projective effective  quasilinear fan is Kähler.
\end{thm}

This theorem produces a large class of Kähler tropical fans. Our companion work~\cite{AP-tht} establishes a hodge theory for Kähler tropical varieties, which are tropical varieties locally modelled by Kähler tropical fans.

\subsection{Tropical Deligne resolution} \label{sec:Deligne-intro}

For a fan $\Sigma$, we denote by $\comp\Sigma$ its canonical compactification, obtained by taking the closure of $\supp\Sigma$ in the tropical toric variety $\TP_\Sigma$.

Let $\Sigma$ be a simplicial homologically smooth tropical fan. For any $\sigma\in \Sigma$, let $\compSigma^\sigma$ be the canonical compactification of the star fan $\ssSigma^\sigma$, and set
\[ H^k(\compSigma^\sigma) \coloneqq \bigoplus_{p+q=k} H^{p,q}(\compSigma^\sigma), \]
for cohomology with rational coefficients. It follows from~\cite[Theorem 1.1]{AP-FY} that $H^k(\compSigma^\sigma,\Q)$ is non-vanishing only in even degrees. Moreover, when $k$ is even, we have $H^k(\compSigma^\sigma,\Q) = H^{k/2,k/2}(\compSigma^\sigma,\Q)$, isomorphic to the piece $A^{k/2}(\ssSigma^\sigma)$ of the Chow ring. In Section~\ref{sec:deligne} we prove the following important theorem.

\begin{thm}[Tropical Deligne resolution]\label{thm:deligne-intro}
  Let $\Sigma$ be a simplicial homologically smooth tropical fan. For each non-negative integer $p$, we have the following long exact sequence
   \[ 0 \rightarrow H^p(\Sigma) \rightarrow \bigoplus_{\sigma \in \Sigma \\
  \dims{\sigma} =p} H^0(\comp \Sigma^\sigma) \rightarrow \bigoplus_{\sigma \in \Sigma \\
  \dims{\sigma} =p-1} H^2(\comp \Sigma^\sigma) \rightarrow \dots \rightarrow \bigoplus_{\sigma \in \Sigma \\
  \dims{\sigma} =1} H^{2p-2}(\comp \Sigma^\sigma) \rightarrow H^{2p} (\comp \Sigma) \to 0. \]
\end{thm}

In the above sequence, the first term $H^p(\Sigma)$ coincides with the coefficient group $\SF^p(\conezero)$ consisting of rational $p$-multiforms on $\Sigma$ (see Section~\ref{sec:prel} for the definition), the notation $\dims\sigma$ stands for the dimension of the face $\sigma$, and the maps between cohomology groups are given by tropical Gysin maps, as we expand in Section~\ref{sec:deligne}. The cohomology groups are with $\Q$-coefficients.

The name given to the theorem above comes from Deligne's weight spectral sequence in Hodge theory~\cite{Del-hodge2} that puts a mixed Hodge structure on the cohomology of a smooth quasi-projective complex variety. In the case of a complement of an arrangement of complex hyperplanes, the cohomology in degree $p$ coincides with the coefficient group $\SF^p(\conezero)$ of the Bergman fan $\ssSigma_\Ma$, for the matroid $\Ma$ associated to the hyperplane arrangement, by the work of Orlik-Solomon~\cite{OT13} and Zharkov~\cite{Zha13}. A wonderful compactification of the hyperplane arrangement's complement produces a long exact sequence as above that involves the cohomology groups of the strata in the compactification. This is a consequence of a theorem of Shapiro~\cite{Sha93} that states that the cohomology of the complement of an arrangement of complex hyperplanes has pure Hodge structure of Hodge-Tate type, of weight $2p$ concentrated in bidegree $(p,p)$ in cohomological degree $p$.

In the theorem above, we obtain a purely polyhedral analogue of the Deligne exact sequence, valid for any homologically smooth tropical fan. In Section~\ref{subsec:Deligne_non_smooth}, we prove a weaker statement for general tropical fans.

\subsection{Reformulation of homological smoothness}

In~\cite{Aks21}, Aksnes proves the following.

\begin{thm}[Aksnes~\cite{Aks21}]
  Let $\Sigma$ be a tropical fan of dimension $d\geq 2$. Suppose that the Borel-Moore homology groups $H_{a,b}^\BM(\ssSigma)$ vanish for all pairs $a,b$ with $b \neq d$. If each of the star-fans $\Sigma^{\sigma}$ for $\sigma \neq \conezero$ satisfies tropical Poincaré duality, then $\Sigma$ satisfies tropical Poincaré duality.
\end{thm}

From the above theorem, proceeding by induction,  we easily deduce the following reformulation of homological smoothness in terms of vanishing of Borel-Moore homology groups.

\begin{thm}
  A tropical fan $\Sigma$ is homologically smooth if and only if it is quasilinear in codimension one, and the vanishing property $H_{a,b}^\BM(\Sigma^\sigma, \Q)=0$ holds for each face $\sigma$ of\/ $\Sigma$ and each pair of integers $a,b$ with $b\neq d- \dims \sigma$.
\end{thm}

Being quasilinear in codimension one means simply that for each face $\sigma$ of codimension one in $\Sigma$, in the one-dimensional tropical fan $\Sigma^\sigma$, there is a \emph{unique} relation between the primitive vectors of rays, given by the balancing condition.

\subsection{Related work}

Tropical homology is introduced in~\cite{IKMZ}. It is notably proved there that in a one-parameter family of smooth projective complex varieties tropicalizing to a tropical variety that is locally modelled by Bergman fans of matroids, the Hodge numbers of smooth fibers are captured by the dimension of tropical cohomology groups of the limit tropical variety. A generalization of this result to more general type of degenerations is obtained in our joint work~\cite{AAPS} wih Aksnes and Shaw. The proof uses the tropical Deligne resolution, Theorem~\ref{thm:deligne-intro}. This theorem plays as well a crucial role in our companion works~\cite{AP-tht, AP20-hc} which establish a Hodge theory for tropical varieties. A similar in spirit resolution dealing with the Stanley-Reisner rings of simplicial complexes and their quotients by (generic) linear forms plays a central role in the recent works~\cites{Adi18, AY20} on partition complexes, and have found several interesting applications. The result is also connected to the recent work~\cite{Pet17} on the cohomology of stratified spaces which provides a generalization of Deligne's spectral sequence. In~\cite{AP-FY} we relate tropical cohomology to combinatorial Hodge theory, we refer to surveys~\cites{Baker18,Huh18, Huh22, Okou22} for a discussion of combinatorial Hodge theory and its applications in combinatorics. Further results in connection to Hodge theory are obtained in~\cite{MZ14, AB14, ARS}.

Finding a classification of tropical fans that verify tropical Poincaré duality is a challenging problem. Interesting results in this direction have been obtained in~\cite{Aks19}. The work~\cite{JSS19} establishes tropical Poincaré duality for tropical varieties locally modelled by Bergman fans of matroids. A sheaf theoretic approach to tropical homology is developed in~\cites{GS-sheaf, GS19}.

\subsection{Organization of the paper}

In Section~\ref{sec:prel} we provide necessary background in polyhedral geometry, and recall the definition of tropical homology and cohomology, tropical modification, and $\mT$-stability. This section contains as well the proof of Theorem~\ref{thm:tropical_homology_manifold_local_fans}. Tropical Deligne resolution for homologically smooth tropical fans is proved in Section~\ref{sec:deligne}.  The proof of Theorem~\ref{thm:smooth_shellable_intro} on $\mT$-stability of homological smoothness will occupy Sections~\ref{sec:coefficients_tropical_modification} and~\ref{sec:homology_tropical_modification}. The technical part involves proving Theorem~\ref{thm:homology_tropical_modification}, that describes the tropical homology of a tropical modification. Final Section~\ref{sec:examples} provides examples related to the results of the paper.

\smallskip
In this article, unless otherwise stated, we will work with homology and cohomology with rational coefficients.

\subsection{Basic notation} \label{sec:intro-basic-notations}

The set of non-negative real numbers is denoted by $\ss\R_{\geq 0}$. In this paper, $N$ will be a free $\Z$-module of finite rank. We denote $M=N^\dual = \hom(N, \Z)$ the dual of $N$. The corresponding rational and real vector spaces are denoted respectively by $\ssN_{\Q}$, $\ssN_{\R}$, and $\ssM_{\Q} = \ssN_{\Q}^\dual$, $\ssM_{\R} = \ssN_{\R}^\dual$.

\smallskip
We use basic results and constructions in multilinear algebra. In particular, for a lattice $N$, we view its dual $M=\ssN^\dual$ as (linear) forms on $N$ and refer to elements of the exterior algebras $\bigwedge^\bul N$ and $\bigwedge^\bul M$ as multivectors and multiforms, respectively. If $\ell \in M$ is a linear map on $N$ and if $\vv \in \bigwedge^p N$ is a multivector, we denote by $\ss\iota_\ell(\vv) \in \bigwedge^{p-1} N$ the contraction of $\vv$ by $\ell$. For $v_1, \dots, v_p \in N$ and $\vv = v_1 \wedge \dots \wedge v_p$, this is given by
\[ \ss\iota_\ell(v_1 \wedge \dots \wedge v_p) = \sum_{i=1}^p (-1)^{i-1} \ell(v_i) v_1 \wedge \dots \wedge \hat{v_i} \wedge \dots \wedge v_p \in \bigwedge^{p-1} N, \]
where the notation $\hat{v_i}$ indicates, as usual, that the factor $v_i$ is removed from the wedge product. We extend this definition to $k$-forms, for any $k\geq2$, by setting, recursively, $\ss\iota_{\ell \wedge \alpha} := \ss\iota_\alpha \circ \ss\iota_\ell$ for any $\alpha \in \bigwedge^k M$ and $\ell\in M$.

\smallskip
We denote by $\T \coloneqq \R \cup \{\ssinfty\}$ the extended real line with the topology induced by that of $\R$ and a basis of open neighborhoods of infinity given by intervals $(a, \ssinfty]$ for $a\in \R$. Extending the addition of $\R$ to $\T$ by setting $\ssinfty + a = \ssinfty$ for all $a \in \T$, endows $\T$ with the structure of a topological monoid called the monoid of tropical numbers. We denote by $\T_+ \coloneqq \ss\R_{\geq 0} \cup\{\ssinfty\}$ the submonoid of non-negative tropical numbers with the induced topology. Both monoids admit a natural scalar multiplication by non-negative real numbers (setting $0\cdot\ssinfty=0$). Moreover, the multiplication map is continuous. As such, $\T$ and $\T_+$ can be seen as modules over the semiring $\ss\R_{\geq 0}$. A polyhedral cone $\sigma$ in $\ssN_{\R}$ is another example of a module over $\R_+$.

\smallskip

The cones appearing in this paper are all strictly convex, meaning that they do not contain any line. If $\sigma$ is a polyhedral cone in $\ssN_{\R}$, we denote by $\ssN_{\sigma, \R}$ the real vector subspace of $\ssN_{\R}$ generated by points of $\sigma$, and define $\ssN_{\R}^\sigma \coloneqq \rquot{\ssN_{\R}}{\ssN_{\sigma, \R}}$. If $\sigma$ is rational, we get natural lattices of full rank $\ssN_\sigma \subset \ssN_{\sigma, \R}$ and $\ssN^\sigma \subset \ssN_{\R}^\sigma$.

Given a poset $(P, \subface)$ and a functor $\phi$ from $P$ to a category $\mathcal C$, if $\phi$ is covariant (\resp contravariant), then for a pair of elements $\tau \subface \sigma$ in $P$, we denote by $\phi_{\tau \subface \sigma}$ (\resp $\phi_{\sigma \supface \tau}$), the corresponding map $\phi(\tau) \to \phi(\sigma)$ (\resp $\phi(\sigma) \to \phi(\tau)$) in $\mathcal C$, the idea being that in the subscript of the map $\phi_{\bul}$ representing the arrow in $\mathcal C$, the first item refers to the source and the second to the target. This convention will be in particular applied to the poset of faces in a fan.

\section{Preliminaries}\label{sec:prel}

\subsection{Fans} \label{subsec:fans}

Consider a fan $\Sigma$ of dimension $d$ in $\ssN_{\R}$. We denote by $\dims\sigma$ the dimension of each cone $\sigma$ in $\Sigma$. We denote by $\ssSigma_k$ the set of $k$-dimensional cones of $\Sigma$. The elements of $\ssSigma_1$ are called \emph{rays}. We denote by $\conezero$ the cone $\{0\}$. A $k$-dimensional cone $\sigma$ in $\Sigma$ is determined by its set of rays. The \emph{support} of $\Sigma$ is denoted by $\supp \Sigma$: it is the closed subset of $\ssN_{\R}$ defined by the union of the cones in $\Sigma$. We call it a \emph{fanfold}. If there is no risk of confusion, we use $\Sigma$ when referring both to the fan and the fanfold. A maximal cone in $\Sigma$ is called a \emph{facet}. We say that $\Sigma$ is \emph{pure dimensional} if all its facets have the same dimension. A fan $\Sigma$ is called \emph{rational} if all its cones are rational. It is called \emph{simplicial} if each cone is generated by as many rays as its dimension.

Given a cone $\sigma\in \Sigma$, the \emph{star fan} $\ssSigma^\sigma$ refers to the fan in $\ssN_{\R}^\sigma=\rquot {\ssN_{\R}}{\ssN_{\sigma,\R}}$ obtained by taking the projection of the cones $\eta$ in $\Sigma$ that contain $\sigma$ as a face, for the projection map $\ssN_{\R} \to \ssN_{\R}^\sigma$. We denote by $\sseta^{\sigma} \in \ssN_{\R}^\sigma$ this projection.

\subsection*{Convention}

We endow the fan $\Sigma$ with the partial order $\subface$ given by the inclusion of cones in $\Sigma$: we write $\tau \subface \sigma$ if $\tau \subseteq \sigma$. We say $\sigma$ \emph{covers} $\tau$ and write $\tau \ssubface \sigma$ if $\tau\subface\sigma$ and $\dims{\tau} = \dims\sigma-1$. The \emph{meet} operation $\wedge$ is defined as follows. For two cones $\sigma$ and $\delta$ of $\Sigma$, we set $\sigma \wedge \delta \coloneqq \sigma \cap \delta$. The set of cones in $\Sigma$ that contain both $\sigma$ and $\delta$ is either empty or has a minimal element $\eta \in \Sigma$. In the latter case, we say that $\eta$ is the \emph{join} of $\sigma$ and $\delta$ and denote it by $\sigma \vee \delta \coloneqq \eta$.

We use $\sigma$ (or any other face of~$\Sigma$) as a superscript when referring to the quotient of some space by $\ssN_{\sigma, \R}$ or to the elements related to this quotient. In contrast, we use $\sigma$ as a subscript for subspaces of $\ssN_{\sigma,\R}$ or for elements associated to these subspaces.

If $\tau\subface \sigma$ are faces of $\Sigma$, we denote by $\sspi_{\tau \subface \sigma}$ all the projection maps $\ssN^\tau\to \ssN^\sigma$, $\ssN^\tau_\Q\to \ssN^\sigma_\Q$, and $\ssN^\tau_{\R}\to \ssN^\sigma_{\R}$.

\subsection{Canonical compactifications}\label{sec:canonical_compactification}

We briefly discuss canonical compactifications of fans. More details can be found in~\cites{OR11, AP-FY}.

Let $\Sigma$ be a fan in $\ssN_{\R}$. For any cone $\sigma$, denote by $\ss\sigma^\vee$ the \emph{dual cone} defined by
\[ \ss\sigma^\vee \coloneqq \left\{m \in \ssM_{\R} \:\st\: \langle m, a \rangle \geq 0 \:\textrm{ for all } a \in \sigma\right\}. \]

The \emph{canonical compactification} $\comp\sigma$ of $\sigma$ is given by $\hom_{\RMod}(\ss\sigma^\vee, \ss\T_{\geq 0})$, \ie, by the set of morphisms $\ss\sigma^\vee \to \ss\T_{\geq 0}$ in the category of $\ss\R_{\geq 0}$-modules. We identify naturally $\sigma$ with the subset $\hom_{\RMod}(\ss\sigma^\vee, \ss\R_{\geq 0})$ of $\comp\sigma$. There is a natural topology on $\comp\sigma$ that makes it a compact topological space. The induced topology on $\sigma$ coincides moreover with the Euclidean topology. We refer sometimes to $\comp\sigma$ as the extended cone of $\sigma$. For a face $\tau \subface \sigma$, we get an inclusion of the compactifications $\comp \tau \subseteq \comp \sigma$ that identifies $\comp \tau$ as the topological closure of $\tau$ in $\comp \sigma$.

The \emph{canonical compactification} $\comp\Sigma$ is defined as the union of $\comp\sigma$, $\sigma\in \Sigma$, where for $\tau \subface \sigma$ in $\Sigma$, we identify $\comp\tau$ with the closure of $\tau$ in $\comp\sigma$. The topology of $\comp\Sigma$ is the induced quotient topology. Each extended cone $\comp\sigma$ naturally embeds as a subspace of $\comp\Sigma$.

There is a special point $\ssinfty_\sigma$ in $\comp\sigma$ defined by the map $\ss\sigma^\vee\to\ss\T_{\geq 0}$ that vanishes on the orthogonal space $\ss\sigma^\perp \coloneqq \left\{m\in M_{\R}\:\st\: \langle m, a \rangle = 0 \:\textrm{ for all } a \in \sigma\right\}$ and takes value $\ssinfty$ everywhere else. Note that for the cone $\conezero$, we have $\ssinfty_{\conezero} = 0$.

The compactification $\comp\Sigma$ naturally lives in the tropical toric variety $\ssTP_\Sigma$ defined as follows. For $\sigma\in \Sigma$, let $\chart\sigma \coloneqq \hom_{\RMod}(\ss\sigma^\vee,\T)$. Since $\hom_{\RMod}(\ss\sigma^\vee, \R)\simeq \ssN_{\R}$, this is a partial compactification of $\ssN_{\R}$. For a pair $\tau \subface \sigma$ in $\Sigma$, we get an inclusion $\~\tau\subseteq\~\sigma$. Gluing the spaces $\chart\sigma$ along these inclusions gives $\ssTP_\Sigma$.

We set $\ssN_{\infty,\R}^\sigma\coloneqq \ssN_{\R} + \ssinfty_\sigma \subseteq \~\sigma$. In this notation, we have $\ssN_{\infty,\R}^\conezero=\ssN_{\R}$. More generally, we have an isomorphism $\ssN_{\infty,\R}^\sigma \simeq \ssN_{\R}^\sigma$.

The tropical toric variety $\ssTP_\Sigma$ is naturally stratified into the disjoint union of \emph{tropical torus orbits} $\ssN_{\infty,\R}^\sigma \simeq \ssN_{\R}^\sigma$, $\sigma \in \Sigma$. The natural inclusion of $\comp\sigma$ into $\chart\sigma$ gives an embedding $\comp\Sigma \subseteq \ssTP_\Sigma$ that identifies $\comp \Sigma$ as the closure of $\Sigma$ in $\ssTP_\Sigma$. We say that a point $x\in \comp\Sigma$ has \emph{sedentarity} $\sigma$ if it belongs to $\ssN_{\infty,\R}^\sigma\cap \comp\Sigma$. In this case, we write $\sed(x)=\sigma$.

The intersection $\comp\Sigma=\ssN_{\infty,\R}^\sigma \simeq \ssN_{\R}^\sigma$ is identified naturally with the star fan $\ssSigma^\sigma$. We denote it by $\ssSigma_\infty^\sigma$. We refer to each face of $\ssSigma_\infty^\sigma$ as a \emph{face of\/ $\comp\Sigma$ of sedentarity $\sigma$.} Note that each face of sedentarity $\sigma$ is in one-to-one correspondance with a face $\eta \supface \sigma$ in $\Sigma$. We denote that face by $\sseta^{\sigma}_\infty$. It is explicitly given by  
\[\sseta^{\sigma}_\infty = \{\ssinfty_{\sigma} + x \mid x\in \eta\} = \comp\eta \cap \ssN^\sigma_{\infty,\R}. \]
We denote by $\cube^\sigma_\eta$ the closure of $\sseta^{\sigma}_\infty$ in $\comp\Sigma$. For each face $\delta=\cube^\sigma_\eta$ of $\comp\Sigma$, we set $\ssN_\delta \coloneqq \ssN_{\infty, \eta}^\sigma \simeq \rquot{\ssN_{\eta}}{\ssN_{\sigma}}$.

\subsection{Unit normal vectors, canonical multivectors, and canonical forms} \label{sec:canonical_forms}

Consider a rational fan $\Sigma$ of dimension $d$. Consider a cone $\sigma$ of $\Sigma$ and let $\tau$ be a face of codimension one in $\sigma$. The hyperplane $\ssN_{\tau,\R} \subset \ssN_{\sigma,\R}$ cuts the space $ \ssN_{\sigma,\R}$ into two closed half-spaces. Denote the half-space that contains $\sigma$ by $\ssH_\sigma$. A \emph{unit normal vector to $\tau$ in $\sigma$} refers to a vector $v$ of $\ssN_\sigma \cap \ssH_\sigma$ that verifies $\ssN_\tau + \Z v = \ssN_\sigma$. We denote by $\nvect_{\sigma/\tau}$ such an element. Note that $\nvect_{\sigma/\tau}$ induces a well-defined generator of $\ssN^\tau_\sigma = \rquot{\ssN_\sigma}{\ssN_\tau}$: we denote this generator by $\sse^\tau_\sigma$. We naturally extend the definition to similar pair of faces in $\comp\Sigma$ having the same sedentarity. In the case $\sigma$ is a ray (and $\tau$ is a point of the same sedentarity as $\sigma$), we also use the notation $\sse_\sigma$ instead of $\nvect_{\sigma/\tau}$.

On each cone $\sigma \in \Sigma$, we fix a generator of $\bigwedge^{\dims{\sigma}}\ssN_\sigma$ that we denote by $\ssnu_\sigma$ and refer to it as the \emph{canonical multivector of $\sigma$} (this is unique up to a sign). Passing to the dual space, we denote by $\ssvarpi_\sigma\in \bigwedge^{\dims{\sigma}}\ssN_\sigma^\dual$ the element that takes value one on $\ssnu_\sigma$. We refer to $\ssvarpi_\sigma$ as the \emph{canonical form of $\sigma$} (again, this is unique up to a sign). We choose $\ssnu_{\conezero} = 1$ and $\ssnu_\rho = \sse_\rho$ for any ray $\rho$.

The above definitions are extended to $\comp\Sigma$ in a natural way by choosing an element $\ssnu^\tau_\sigma \in \bigwedge^{\dims{\sigma}-\dims{\tau}}\ssN^\tau_{\infty,\sigma}\simeq \bigwedge^{\dims{\sigma}-\dims{\tau}}\ssN^\tau_\sigma$ for each $\cube^\tau_\sigma$, $\tau \subface \sigma$, and setting $\varpi^\tau_\sigma \coloneqq \ssnu^{\tau\,\dual}_\sigma$. In the case of interest to us, $\Sigma$ is simplicial and in this case, there is a natural choice for the extension: if $\tau \subface \sigma$ is a pair of cones in $\Sigma$, there exists a unique minimal $\tau'$ such that $\sigma = \tau \vee \tau'$. We choose $\ssnu^\tau_\sigma$ to be the image of $\ssnu_{\tau'}$ under the projection map from $N$ to $\ssN_{\infty, \R}^\tau\simeq \ssN^\tau$.

\subsection{Orientation and sign function} \label{subsec:orientation}

Using the choices of canonical multivectors discussed in the previous section, we give an orientation (in the topological sense) to each face of $\comp\Sigma$. Using this, we associate a \emph{sign} denoted by $\sign(\gamma,\delta)$ to each pair of closed faces $\gamma \ssubface \delta$ in $\comp\Sigma$,  as follows.

In the case both the faces $\gamma$ and $\delta$ have the same sedentarity, we define $\sign(\gamma,\delta)$ to be the sign of the evaluation $\ssvarpi_\delta\bigl(\nvect_{\delta/\gamma} \wedge \ssnu_\gamma\bigr)$. Otherwise, there will exist a pair of cones $\tau' \ssubface \tau$ of $\Sigma$, and $\sigma \in \Sigma$ such that $\gamma = \cube^{\tau'}_\sigma$ and $\delta = \cube^{\tau}_\sigma$. Consider the projection map $\ss\sspi_{\gamma \ssubface \delta} \colon \ssN_\gamma = \ssN_{\infty, \sigma}^{\tau'} \to \ssN_{\delta}=\ssN_{\infty, \sigma}^\tau$ between the two corresponding lattices of $\gamma$ and $\delta$. This map is surjective, and induces a surjective linear map $\bigwedge^{\bul}\ssN_\gamma \to \bigwedge^{\bul}\ssN_{\delta}$. There exists an element $\ssnu'$ such that $\sspi_{\gamma\ssubface\delta} (\ssnu') = \ssnu_{\gamma}$. Consider the primitive vector $\sse^{\tau'}_\tau$ of the ray $\sstau^{\tau'}_{\infty}$ in $\sssigma^{\tau'}_{\infty}$. Then, $\sign(\gamma,\delta)$ is defined to be the sign of $-\ssvarpi_\delta(\sse^{\tau'}_\tau \wedge \ssnu')$. Although we made a choice for $\ssnu'$, the element $\sse^{\tau'}_\tau \wedge \ssnu'$ will be independent of that choice. In other words, the sign function is well-defined.

\subsection{Tropical homology and cohomology groups} \label{subsec:homology}
Let $\Sigma$ be a rational fan and $\comp\Sigma$ its canonical compactification. The (extended) polyhedral structures on $\Sigma$ and $\comp\Sigma$ allow to associate tropical homology groups and cohomology groups to them, as introduced in the work \cite{IKMZ}. We refer to~\cite{JSS19, AP-tht} for more details.

We first recall the definition of the multi-tangent spaces $\SF_p$, and their duals ~$\SF^p$ called multi-cotangent spaces, associated to faces, working with rational coefficients. We do this for $\comp\Sigma$ with the face structure given by $\cube^\tau_\sigma$, for pairs of faces $\tau \subface \sigma$ of $\Sigma$. This provides cellular chain and cochain complexes that calculate tropical homology and tropical cohomology groups of the canonical compactification $\comp \Sigma$. It is easy to adapt the definitions to $\Sigma$ with the natural face structure given by its cones.

For each face $\gamma = \cube^\tau_\sigma$ of $\comp\Sigma$, we denote by $\ssN_{\gamma,\Q}$ the rational vector space defined by the lattice $\ssN_{\gamma}$. For each integer $p\in\ss\Z_{\geq 0}$, we define the \emph{$p$-th multi-tangent} and the \emph{$p$-th multi-cotangent spaces of\/ $\comp\Sigma$ at $\gamma$}, with rational coefficients, denoted by $\SF_p(\gamma)$ and $\SF^p(\gamma)$ respectively, by
\[ \SF_p(\gamma)\coloneqq \hspace{-.3cm} \sum_{\delta \supface \gamma \\ \sed(\delta) = \tau } \hspace{-.3cm} \bigwedge^p \ssN_{\delta,\Q} \ \subseteq \bigwedge^p \ssN^\tau_\Q, \quad \textrm{and} \quad \SF^p(\gamma) \coloneqq \SF_p(\gamma)^\dual. \]

An inclusion of faces $\gamma \subface \delta$ in $\comp\Sigma$, induces maps $\ssi_{\delta \supface \gamma}\colon \SF_p(\delta) \to \SF_p(\gamma)$ and $\ssi_{\gamma \subface \delta}^*\colon \SF^p(\gamma) \to \SF^p(\delta)$. These are defined as follows. In the case $\gamma$ and $\delta$ have the same sedentarity, then $\ssi_{\delta \supface \gamma}$ is just an inclusion. In the case $\gamma=\cube_\sigma^{\tau'}$ and $\delta=\cube_\sigma^\tau$ for cones $\tau \subface \tau' \subface \sigma$ in $\Sigma$, the map $\ssi_{\delta\supface \gamma}$ is induced from the projection map $N^\tau_{\infty} \to N^{\tau'}_{\infty}$. In the general case, $\ssi_{\delta\supface \gamma}$ is the composition of a projection and an inclusion; the map $\ssi_{\gamma \subface \delta}^*$ is the dual of $\ssi_{\delta\supface \gamma}$.

Suppose tha $X = \Sigma$ or $\comp\Sigma$. For a pair of  $p,q \in \ss\Z_{\geq 0}$, we define 
\[ C_{p,q}(X) \coloneqq \bigoplus_{\gamma \in X \\ \dims\gamma =q \\ \gamma \text{ compact}} \SF_p(\gamma). \]
For each non-negative integer $p$, we get a chain complex 
\[ C_{p,\bullet}(X)\colon \quad \dots\longrightarrow C_{p, q+1}(X) \xrightarrow{\ \partial_{q+1}\ } C_{p,q}(X) \xrightarrow{\ \partial_q\ } C_{p,q-1} (X)\longrightarrow\cdots \]
with the differential given by the sum of maps $\sign(\gamma,\delta)\cdot\ssi_{\delta\ssupface \gamma}$, with the sign function as in Section~\ref{subsec:orientation}.

\smallskip
The tropical homology of $X$ with rational coefficients  is defined by
\[ H_{p,q}(X) \coloneqq H_{p,q}(X,\Q)\coloneqq H_q(C_{p,\bullet}(X)). \]
Proceeding similarly, for each non-negative integer $p$, we get a cochain complex
\[ C^{p,\bullet}(X)\colon \quad \dots\longrightarrow C^{p, q-1}(X) \xrightarrow{\ \d^{q-1}\ } C^{p,q}(X) \xrightarrow{\ \d^{q}\ } C^{p,q+1}(X) \longrightarrow \cdots \]
with $(p,q)$ term given
\[ C^{p,q}(X) \coloneqq C_{p,q}(X)^\dual \simeq \! \bigoplus_{\gamma \in X \\ \dims\gamma = q \\ \gamma \text{ compact}} \!\! \SF^p(\delta). \]
The tropical cohomology of $X$  with rational coefficients is defined by
\[ H^{p,q}(X) \coloneqq H^{p,q}(X, \Q)\coloneqq H^q(C^{p,\bullet}(X)). \]

The compact-dual versions of tropical homology and cohomology are defined by allowing non-compact faces. These are called Borel-Moore homology and cohomology with compact support, defined as follows. First, we define
\[ C^\BM_{p,q}(X) \coloneqq \bigoplus_{\gamma \in X \\ \dims\gamma=q} \SF_p(\gamma) \quad \text{and} \quad C_c^{p,q}(X) \coloneqq \bigoplus_{\gamma \in X \\ \dims\gamma = q } \SF^p(\gamma). \]
Together, they give (co)chain complexes $C^\BM_{p,\bul}(X)$ and $C_c^{p,\bul}(X)$. The (co)homology of these complexes are called the tropical Borel-Moore homology and the tropical cohomology with compact support with rational coefficients, respectively,
\[ H^\BM_{p,q}(X) \coloneqq H_q(C^\BM_{p,\bul}(X)) \quad \text{and} \quad H_c^{p,q}(X) \coloneqq H^q(C_c^{p,\bul}(X)). \]
Both notions of homology and both notions of cohomology coincide in the case $X$ is compact. 

We define relative versions of the chain and cochain complexes, and their homology and cohomology, respectively, for a pair $\Delta \subset \Sigma$ of fans. For example, the relative complex $C^\BM_{p,\bul}(\Sigma, \Delta)$ is defined as the cokernel of the inclusion $C^\BM_{p,\bul}(\Delta) \hookrightarrow C^\BM_{p,\bul}(\Sigma)$, and
\[H^\BM_{p,q}(\Sigma, \Delta) \coloneqq H_q(C^\BM_{p,\bul}(\Sigma, \Delta)).\]

\smallskip

In this paper, homology and cohomology refer to the tropical ones, and for the sake of simplicity, we omit  the word tropical. We note that as in the classical setting, there exist other ways of computing the same groups. In particular, a sheaf theoretic version of these constructions allows to associate cohomology and homology groups to open subsets of $\comp\Sigma$, and more general tropical varieties, we refer to~\cite{JSS19,GS-sheaf} for a discussion. Homology and cohomology depends only on the support.

\smallskip
The collection of coefficient sheaves and cosheaves come with contraction maps defined as follows.

Let $\sigma \in \Sigma$ and consider an element $\ssnu \in \SF_{k}(\sigma)$. Given an element $\alpha \in \SF^{p}(\sigma)$, we denote by $\ss\contract_{\ssnu}(\alpha)$ the element of $\SF^{p-k}(\sigma)$ that on each element $\ssnu' \in \SF_{p-k}(\sigma)$ takes value $\alpha(\ssnu \wedge \ssnu')$. The map $\ss\contract_{\ssnu} \colon \SF^{p}(\sigma) \to \SF^{p-k}(\sigma)$ is linear. We extend the definition naturally to all faces $\delta$ of $\comp\Sigma$ and obtain maps $\ss\contract_{\ssnu} \colon \SF^{p}(\delta) \to \SF^{p-k}(\delta)$ for any element $\ssnu \in \SF_k(\delta)$. We call $\ss\contract_{\ssnu}(\alpha)$ the contraction of $\alpha$ by $\ssnu$, and refer to $\ss\contract_{\ssnu}$ as the contraction map defined by $\ssnu$. Dually, for $\alpha \in \SF^{k}(\sigma)$, we get contraction maps $\contract_\alpha \colon \SF_p(\sigma) \to \SF_{p-k}(\sigma)$.

These maps allow to define a cap product between homology and cohomology groups.

\subsection{Tropical Feichtner-Yuzvinsky}

Let $\Sigma$ be a rational simplicial fan in $\ssN_{\R}$, and denote by $\comp\Sigma$ its canonical compactification. Denote by $A^\bul(\Sigma)$ the Chow ring of $\Sigma$ with rational coefficients, defined as the Chow ring of the toric variety associated to $\Sigma$. The following result provides a description of the cohomology of the canonical compactification $\comp\Sigma$.

\begin{thm}[Tropical Feichtner-Yuzvinsky for fans~\cite{AP-FY}] \label{thm:ring_morphism_homology}
For any integer $p$, there is an isomorphism
\[ \begin{array}{ccc}
 H^{p,p}(\comp\Sigma) & \longsimto & A^p(\Sigma).
\end{array} \]
Altogether, they induce a ring morphism $H^{\bul,\bul}(\comp\Sigma) \to A^\bul(\Sigma)$, by mapping $H^{p,q}(\comp\Sigma)$ to zero in the bidegree $p\neq q$. In addition,  $H^{p,q}(\comp\Sigma)$ is trivial for $p < q$ and for $p>q=0$.
\end{thm}

\subsection{Tropical fans} \label{sec:tropical-fans}

A \emph{tropical fan} in $N$ is a pair $(\Sigma, \ssomega_\Sigma)$ consisting of a pure dimensional rational fan $\Sigma$ of dimension $d$, and a tropical orientation $\ssomega_\Sigma$, that is,
an integer valued map
\[ \ssomega_\Sigma\colon \Sigma_d \to \Z\setminus\{0\} \]
that verifies
\[ \sum_{\sigma \ssupface \tau} \ssomega_\Sigma(\sigma)\sse_{\sigma}^{\tau} = 0 \qquad \textrm{in} \qquad \ssN^\tau \]
for each cone $\tau$ in $\Sigma$ of codimension one.
We refer to the support of $(\Sigma, \ssomega_\Sigma)$ as a \emph{tropical fanfold} and denote it by $\supp{\Sigma}$ or simply $\Sigma$ if there is no risk of confusion. We endow $\supp{\Sigma}$ with the tropical orientation induced by $\ssomega_\Sigma$ on the set of regular points (those points having an open neighborhood isomorphic to an open subset of a real vector space).

For any face $\sigma \in \Sigma$, the star fan $\ssSigma^{\sigma}$ gets the induced orientation $\ssomega_{\Sigma^\sigma}$ defined by $\ssomega_{\Sigma^\sigma}(\sseta^{\sigma}) = \ssomega_\Sigma(\eta)$, for any facet $\eta \supface \sigma$. The star fan at $\sigma$ of a tropical fan $(\Sigma, \ssomega_\Sigma)$ refers to $(\ssSigma^\sigma, \ssomega_{\Sigma^\sigma})$

The balancing condition ensures the existence of a \emph{fundamental class} $[\Sigma]\in H_{d,d}^\BM(\Sigma)$. This is defined by using canonical multivectors (see Section~\ref{sec:canonical_forms}). Concretely, the class $[\Sigma]$ is represented by \emph{the canonical element} $\ssnu_{\Sigma}$ of $C_{d,d}^{\BM}(\Sigma)$ given by 
\[ \ssnu_{\Sigma} \coloneqq \bigl(\ssomega_\Sigma(\sigma)\ssnu_\sigma\bigr)_{\sigma \in \Sigma_d} \in \bigoplus_{\sigma \in \Sigma_d} \bigwedge^d \ssN_{\sigma,\Q}. \]
(Note that the canonical element and canonical class are both integral.)

Taking the cap product with the fundamental class gives a map
\[ \frown [\Sigma] \colon H^{p,q}(\Sigma)\to H^{\BM}_{d-p,d-q}(\Sigma), \]
for each pair of integers $p,q\in \{0,\dots,d\}$. A tropical fan $\Sigma$ satisfies \emph{tropical Poincaré duality} if these maps are all isomorphisms.

The above maps are always injective for any tropical fan. Moreover, $H^{p,q}(\Sigma)$ is trivial for $q \neq 0$, and for $q=0$, we have $H^{p,q}(\Sigma) =\SF^p(\conezero)$. Tropical Poincaré duality is thus equivalent to the vanishing of $H_{a, b}^\BM(\Sigma)$ for $b\neq d$, and the surjectivity of the natural maps
\[ \SF^p(\conezero) \to H^\BM_{d-p,d}(\Sigma), \qquad \alpha \mapsto \iota_\alpha(\ssnu_{\Sigma}). \]

Cap product with fundamental class is described using contraction maps defined in the previous section. In the case $q=0$, it is the map
\begin{equation}\label{eq:reformulation_PD}
  \SF^k(\conezero) \to H^\BM_{d-k, d}(\Sigma)
\end{equation}
defined as follows. Consider the fundamental class $[\Sigma] \in H^\BM_{d,d}(\Sigma)$. The coefficient of a facet $\sigma\in \Sigma_d$ in $[\Sigma]$ is given by $\omega(\sigma)\ssnu_\sigma$. Given an element $\alpha \in \SF^k(\conezero)$, the corresponding element in $H^\BM_{d-k, d}(\Sigma)$ has coefficient at the facet $\sigma \in \Sigma_d$ given by $\contract_{\alpha}(\omega(\sigma)\ssnu_{\sigma})$. Tropical Poincaré duality for $\Sigma$ requires that these maps being all isomorphisms. Note that, dually, we get a map
\begin{equation}\label{eq:reformulation_PD_dual}
  H_c^{d-k, d}(\Sigma) \to \SF_k(\conezero).
\end{equation}
Since the dual map given in~\eqref{eq:reformulation_PD} is injective for any tropical fan, this map is always surjective. In the case $\Sigma$ verifies tropical Poincaré duality, this is an isomorphism.

A tropical fan $\Sigma$ is called \emph{homologically smooth} if any open subset $U$ of $\supp{\Sigma}$ satisfies tropical Poincaré duality. Equivalently, as claimed in Theorem~\ref{thm:tropical_homology_manifold_local_fans}, a tropical fan $\supp{\Sigma}$ is homologically smooth if and only if each star fan $\ssSigma^\sigma$ of $\Sigma$, $\sigma\in \Sigma$, verifies tropical Poincaré duality.

\begin{proof}[Proof of Theorem~\ref{thm:tropical_homology_manifold_local_fans}]
  Suppose that $\Sigma$ is homologically smooth. Let $\sigma$ be a face of $\Sigma$. A point in the relative interior of $\sigma$ has an open neighborhood $U$ isomorphic to $W \times \ssSigma^\sigma$ for an open neighborhood $W$ of the origin in $\R^{\dims\sigma}$. Since the coefficient sheaves in $\R^{\dims\sigma}$ are constant, it is easy to show that $W$ verifies tropical Poincaré duality. By homological smoothness of $\Sigma$, $U$ verifies tropical Poincaré duality. Applying the Künneth formula, we deduce that $\ssSigma^\sigma$ verifies Poincaré duality. This implies one direction.

  The other direction follows from the arguments given in~\cite[Section 4]{JSS19} for the proof of Theorem 4.33 in \emph{loc.cit.}, replacing \emph{matroidal} condition on star fans by requiring them to verify \emph{tropical Poincaré duality}. We omit the details.
\end{proof}

The following proposition is a consequence of  Theorem~\ref{thm:ring_morphism_homology}, and the compatibility between the duality pairing in the tropical cohomology and the duality pairing in the Chow ring.
\begin{prop}\label{prop:FY_duality} The following assertions are equivalent for a simplicial tropical fan $\Sigma$:
  \begin{itemize}
   \item $\comp\Sigma$ verifies Poincaré duality.
   \item the Chow ring of\/ $\Sigma$ verifies Poincaré duality, and all the cohomology groups $H^{p,q}(\comp\Sigma)$ for $p>q$ vanish.
  \end{itemize}
\end{prop}

Homological smoothness is closed under products. In fact, the following stronger property holds. For two tropical fanfolds $\Sigma_1, \Sigma_2$, $\Sigma_1 \times \Sigma_2$ is homologically smooth if and only if $\Sigma_1$ and $\Sigma_2$ both are.

\subsection{Gysin maps}\label{sec:Gysin}

Consider a simplicial tropical fan $\Sigma$. Assume that $\Sigma$ is homologically smooth. For a cone $\sigma \in \Sigma$, the set of points of sedentarity $\sigma$ form a fan in $\ssN_{\infty,\R}^\sigma$ denoted by $\ssSigma_\infty^\sigma$. This is based at the point $\ssinfty_\sigma$ of $\compSigma$ and is naturally isomorphic to $\ssSigma^\sigma$. Its closure $\compSigma_\infty^\sigma$ in $\compSigma$ is identified naturally with the canonical compactification $\compSigma^\sigma$ of $\ssSigma^\sigma$. This means that the canonical compactification $\compSigma^\sigma \simeq \compSigma_\infty^\sigma$ naturally lives in $\compSigma$. Moreover, for an inclusion of faces $\tau \ssubface \sigma$ of $\Sigma$, we get an inclusion $ \i!_{\sigma \ssupface \tau}\colon \compSigma_\infty^\sigma \hookrightarrow \compSigma_\infty^\tau$. This induces the map on cohomology $\i!_{\tau \ssubface \sigma}^* \colon H^\bul(\compSigma_\infty^\sigma) \to H^\bul(\compSigma_\infty^\tau)$. The tropical varieties $\comp \Sigma_\infty^\sigma$ and $\comp \Sigma_\infty^\tau$ are both tropical homology manifolds. Therefore, by Theorem~\ref{thm:PD-intro}, they both verify tropical Poincaré duality. Taking the dual of $\i!^*_{\sigma \ssupface \tau}$ and using Poincaré duality, we obtain the \emph{Gysin map}
\[ \gys!_{\sigma \ssupface \tau}\colon H^\bul(\comp \Sigma_\infty^\sigma) \to H^{\bul+2}(\comp \Sigma_\infty^\tau). \]

\subsection{Tropical modification} \label{subsec:tropical_modification}

We briefly recall the definition of tropical modification, introduced in the pioneering work of Mikhalkin~\cites{Mik06, Mik07}. We refer to~\cites{BIMS,Kal15}, and~\cite[Section 5]{AP-hodge-fan} for a more through discussion.

Let $f$ be a meromorphic function on a tropical fan $\Sigma$. This means $f$ is a continuous conewise integral linear function on $\Sigma$. The divisor $D=\div(f)$ is a Minkowski weight of dimension $d-1$ defined by the orders of vanishing $\ord_\tau(f)$ of $f$ at $(d-1)$ dimensional cones $\tau$ of $\Sigma$, see~\cite[Section 3]{AR10} and \cite[Section 4]{AP-FY}. Here,
\[
  \ord_\tau(f) \coloneqq -\sum_{\sigma \ssupface \tau} \ssomega_\Sigma(\sigma)\ssf_{\sigma}(\nvect_{\sigma/\tau}) + \ssf_\tau\Bigl(\sum_{\sigma \ssupface \tau} \ssomega_\Sigma(\sigma) \nvect_{\sigma/\tau}\Bigr),
\]
where $\ssf_{\sigma}$ is an integral linear function with $\ssf_\sigma\rest{\sigma} = f\rest{\sigma}$.

Denote by $\Delta$ the subfan of $\Sigma$ given by the support of $\ord(f)$, endowed with tropical orientation $\ssomega_\Delta = \ord(f)$. We allow the case where the divisor $\div(f)$ is trivial, in which case, $\Delta$ will be empty.

The tropical modification of $\Sigma$ induced by $f$ will be a fan in $\~N_\R \simeq \R^{n+1}=\R^{n} \times \R$ that we will denote by $\tropmod{f}{\Sigma}$.  To define it, we first consider the graph $\Gamma!_f$ of $f$ defined by
\[ \begin{array}{rccc}
\Gamma!_f\colon & \supp\Sigma & \longrightarrow & \~N_\R=N_\R \times \R, \\
                & x           & \longmapsto     & (x, f(x)).
\end{array} \]
For $\sigma \in \Sigma$, set $\basetm\sigma\coloneqq \Gamma!_f(\sigma)$, the image in $\~N_\R$ of $\sigma$ by $\Gamma!_f$. For $\delta \in \Delta$, set $\uptm\delta \coloneqq \basetm\delta + \ss\R_{\geq 0} \etm$ where $\etm = (0, 1)\in N_\R \times \R$ ($0$ in $\etm=(0,1)$ is the origin in $\ssN_\R$).

We define the \emph{tropical modification $\tropmod{f}{\Sigma}$ of\/ $\Sigma$ along $\Delta$ with respect to $f$}, simply \emph{tropical modification of $\Sigma$ along $\Delta$} if everything is understood from the context, to be the fan in $\~N_\R\simeq \R^{n+1}$ given by the collection of cones 
\[ \tropmod{f}{\Sigma} \coloneqq \bigl\{\,\basetm\sigma \mid \sigma\in\Sigma\,\bigr\} \cup \bigl\{\,\uptm\delta \mid \delta\in\Delta\,\bigr\}.
\]
Note that $\dims{\uptm\delta} =\dims\delta+1$ for $\delta\in \Delta$, and $\dims{\basetm \sigma} = \dims\sigma$ for $\sigma \in \Sigma$, so that setting $\~\Sigma\coloneqq \tropmod{f}{\Sigma}$, we have
\[ \ss{\~\Sigma}_k = \bigl\{\basetm\sigma \mid \sigma\in\Sigma_k\bigr\} \cup \bigl\{\,\uptm\delta \mid \delta\in\Delta_{k-1}\,\bigr\}. \]

The fan $\~\Sigma$ is endowed with the map $\~\omega \colon \ss{\~\Sigma}_d \to \Z \setminus \{0\}$ defined by $\~\omega(\basetm\sigma) \coloneqq\omega!_\Sigma(\sigma)$ for $\sigma \in \Sigma_d$, and $\~\omega(\uptm\delta) \coloneqq\omega!_\Delta(\delta)$ for $\delta\in \Delta_{d-1}$. It is easy to see that the map $\~\omega$ is a tropical orientation of\/ $\~\Sigma$. This means $\tropmod{f}{\Sigma}$ endowed with $\~\omega$ is a tropical fan. We moreover have a natural conewise integral linear projection
\[ \prtm \colon \supp{\tropmod{f}{\Sigma}} \to \supp\Sigma. \]

We will need the following characterization of the star fans in the tropical modification. Given a cone $\sigma\in \Sigma$, let $\ell \in M$ be an integral linear map which coincides with $f$ on $\sigma$. The difference $f-\ell$ vanishes on $\sigma$. Denote by $\sspi^\sigma\colon N \to \ssN^\sigma$ the projection map. Then, $f-\ell$ induces a meromorphic function on $\ssSigma^\sigma$ denoted by $\ssf^\sigma = \sspi^\sigma_*(f - \ell)$. We refer to $f^\sigma$ as the \emph{meromorphic function induced by $f$ on $\ssSigma^\sigma$}, although it is only defined up to addition by an integral linear function on $\ssN^\sigma$.

\begin{prop} \label{prop:star_fans_tropical_modifications}
Let $\Sigma$ be a tropical fan. Consider a meromorphic function $f$ on $\Sigma$, and denote by $\Delta$ the tropical fan given by $\div(f)$. Let $\~\Sigma \coloneqq \tropmod{f}{\Sigma}$. Then, we have the following description of the star fans in\/ $\~\Sigma$:
\begin{itemize}
  \item If $\delta \in \Delta$, then $\ss{\~\Sigma}^{\uptm\delta} \simeq \ss\Delta^\delta$.
  \item If $\sigma \in \Sigma$, then we have $\ss{\~\Sigma}^{\basetm\sigma} \simeq \tropmod{\ssf^\sigma}{\ssSigma^\sigma}$, where $f^\sigma$ is the meromorphic function induced by $f$ on $\Sigma^\sigma$.
\end{itemize}
\end{prop}

The proof is obtained by a direction verification and is omitted.

\subsection{$\mT$-stability} \label{subsec:T_stability}

In our work~\cite{AP-hodge-fan}, we introduce a notion of \emph{$\mT$-stability} for tropical fans and their geometric properties. This allows in practice to produce new tropical fans with nice properties out of the existing ones. We give a brief presentation here and refer to~\cite{AP-hodge-fan} for more details.

We have three types of operations on tropical fans: products, stellar subdivisions and their inverse stellar assemblies (see~\cite[Section 2.5]{AP-hodge-fan} for the definition), and tropical modifications. The first two are classical, the third, explained in the previous section, is specific to the tropical setting.

Consider a  class of tropical fans $\Csh$ and let $\Sh\subseteq \Csh$ be a subclass. $\Sh$ is called \emph{$\mT$-stable in $\Csh$}, simply \emph{$\mT$-stable} in the case $\Csh$ is the class of all tropical fans, if the following statements are verified:
\begin{itemize}
  \item (Stability under products) For a pair of elements $(\Sigma, \ssomega_\Sigma)$ and $(\Sigma', \ssomega_{\Sigma'})$ in $\Sh$, if $(\Sigma \times \Sigma', \omega!_{\Sigma\times \Sigma'})$ belongs to $\Csh$, then it lies in $\Sh$.

  \item (Stability under a tropical modification along a divisor in the subclass) Given an element $(\Sigma, \omega!_\Sigma)$ in $\Sh$, and a meromorphic function $f$ on $\Sigma$ with $\div(f)$ in $\Sh$, the tropical modification $(\~\Sigma, \omega!_{\~\Sigma})$ of $\Sigma$ with respect to $f$ is in $\Sh$ provided that it belongs to $\Csh$.

  \item (Stability under stellar subdivisions and  assemblies with center in the subclass) For a tropical fan $(\Sigma, \omega!_\Sigma)$ in $\Csh$, and for $\sigma \in \Sigma$ with the property that the stellar subdivision $\Sigma'$ of $\Sigma$ at $\sigma$ belongs to $\Csh$ and $\Sigma^\sigma$ belongs to $\Sh$, we have 
  \[(\Sigma, \omega!_{\Sigma}) \in \Sh \quad \textrm{if and only if} \quad (\Sigma', \omega!_{\Sigma'}) \in \Sh.\]
\end{itemize}

Examples of classes $\Csh$ of interest in this paper are \emph{all}, \resp \emph{simplicial}, \resp \emph{quasi-projective}, and  \resp \emph{div-faithful}  tropical fans, see~\cite{AP-hodge-fan} for a more exhaustive list. We recall from \emph{loc.\;cit.}\ that a tropical fan $\Sigma$ is div-faithful if it verifies the following property: for each $\sigma \in \Sigma$, if the divisor of a meromorphic function defined on the star fan $\ssSigma^\sigma$ is zero, then that meromorphic function is linear.

The subclass $\Csh$ of $\Csh$ is trivially $\mT$-stable in $\Csh$. Moreover, it is easy to see that an intersection of $\mT$-stable classes is again $\mT$-stable. This motivated the following definition.

Consider a subset $\Bsh$ of $\Csh$, viewed as the base set. The \emph{$\mT$-stable subclass of\/ $\Csh$ generated by $\Bsh$} denoted by $\gst{\Bsh}{\Csh}$ is the smallest $\mT$-stable subclass of $\Csh$ that contains $\Bsh$; it is simply obtained by taking the intersection of all subclasses $\Sh$ of $\Csh$ that contain $\Bsh$ and which are $\mT$-stable. In the case $\Csh$ is the class of all tropical fans, we simply write $\gst{\Bsh}{}$.

An important base set is the set $\Bsho$ defined as follows. Denote by $\conezero$ the cone~$\{0\}$. Denote, by an abuse of the notation, by $\conezero$ the fan consisting of unique cone $\conezero$. Let $(\conezero, n)$, $n\in\Z\setminus\{0\}$, be the tropical fan of dimension $0$ with weight function taking value $n$ on $\conezero$. Denote by  $\Lambda$ the complete tropical fan in $\R$ that consists of three cones $\conezero, \R_{\geq 0}$, and $\R_{\leq 0}$, with tropical orientation given by the constant function $1$ on $\R_{\geq 0}$ and $\R_{\leq 0}$. Set
\[ \Bsho\coloneqq \left\{ (\conezero,n) \st  n\in\Z\setminus\{0\} \right\} \cup \left\{ \Lambda \right\}. \]
In this very simple case already, $\gst{\Bsho}{}$ contains many interesting fans. For instance,  Bergman fans of matroids, and generalized Bergman fans (those with support a Bergman fanfold) are all in this class, but $\gst{\Bsho}{}$ is strictly larger.

\smallskip
A tropical fan is called \emph{quasilinear} if it belongs to $\gst{\Bsho}{}$. This means that the tropical fan can be obtained from the collection $\Bsho$ by performing a sequence of the three above operations on tropical fans. More precisely, a tropical fan is quasilinear if it either belongs to~$\Bsho$, or, is a product of two quasilinear fans, or, is the stellar subdivision or stellar assembly of a quasilinear fan, or, is obtained as a result of tropical modification of a quasilinear fan along a tropical divisor which is itself quasilinear.

\smallskip
$\mT$-stability can be defined for properties of tropical fans. If $P$ is a predicate on tropical fans and $\Csh$ is a class of tropical fans, then $P$ is called \emph{$\mT$-stable in $\Csh$} if the subclass of tropical fans in $\Csh$ which verify $P$ is $\mT$-stable in $\Csh$.

\section{Tropical Deligne resolution} \label{sec:deligne}

Let $(\Sigma, \omega!_\Sigma)$ be a simplicial homologically smooth tropical fan. For each cone $\sigma\in \Sigma$, we denote by $\ssSigma^{\sigma}$ the star fan of $\sigma$ in $\Sigma$ that we endow with the induced orientation $\ssomega_{\Sigma^{\sigma}}$. This is a homologically smooth tropical fan in $\ssN^\sigma_\R$. We denote by $\compSigma^{\sigma}$ its canonical compactification. For the canonical compactification of a tropical fan $\Sigma$, recall as well that we set
\[ H^k(\comp \Sigma) := \bigoplus_{p+q=k} H^{p,q}(\comp \Sigma). \]
By Theorem~\ref{thm:PD-intro}, $\comp\Sigma$ verifies tropical Poincaré duality. It follows then from Theorem~\ref{thm:ring_morphism_homology} and Proposition~\ref{prop:FY_duality} that $H^k(\comp \Sigma)$ is non-vanishing only in even degrees, and in that case, it is equal to $H^{k/2,k/2}(\comp \Sigma)\simeq A^{k/2}(\Sigma)$, where $A^\bullet(\Sigma)$ refers to the Chow ring of $\Sigma$ with rational coefficients.

\smallskip
The aim of this section is to prove Theorem~\ref{thm:deligne-intro}, namely that, we have the following long exact sequence
\[ 0 \rightarrow \SF^p(\conezero) \rightarrow \bigoplus_{\sigma \in \Sigma_p} H^0(\compSigma_\infty^\sigma) \rightarrow \bigoplus_{\sigma \in \Sigma_{p-1}} H^2(\compSigma_\infty^\sigma) \rightarrow \dots \rightarrow \bigoplus_{\sigma \in \Sigma_1} H^{2p-2}(\compSigma_\infty^\sigma) \rightarrow H^{2p} (\compSigma) \to 0. \]
The maps between cohomology groups in the above sequence are given by the sum of Gysin maps, see Section~\ref{sec:Gysin}.

Using tropical Poincaré duality for canonical compactifications $\comp \Sigma^\sigma$, consequence of Theorem~\ref{thm:PD-intro}, it will be enough to prove the exactness of the following complex for each $k$ (here $k=d-p$):
\begin{equation} \label{eqn:Deligne_dual}
0 \rightarrow H^{2k}(\compSigma) \to \bigoplus_{\sigma \in \Sigma_1} H^{2k}(\compSigma^\sigma_\infty) \to \bigoplus_{\sigma \in \Sigma_2} H^{2k}(\compSigma^\sigma_\infty) \to \dots \to \bigoplus_{\sigma \in \Sigma_{d-k}} H^{2k}(\compSigma^\sigma_\infty) \to \SF_{d-k}(\conezero) \to 0.
\end{equation}

In order to simplify the notation, we drop sometimes $\infty$ from indices and identify $\compSigma^\sigma$ with $\compSigma_\infty^\sigma$ living naturally in $\comp \Sigma$ as the closure of those points that have sedentarity equal to $\sigma$.

In the next three sections, we give a proof of this theorem using Proposition~\ref{prop:exactness_Omega} that gives a resolution of the sheaf of tropical holomorphic forms, a result of independent interest. We sketch a cellular version of this resolution in Section~\ref{sec:proof_cellular}, that can be used in the same way to get the theorem. In Section \ref{subsec:Deligne_non_smooth}, we prove a weaker statement for any tropical fan.

\subsection{The sheaf of tropical holomorphic forms} \label{sec:proof_deligne_sheaf}

Let $Z = \Sigma$ or $\comp \Sigma$, for a simplicial tropical fan $\Sigma$. We denote by $\ss\Omega^k_Z$ the sheaf of \emph{tropical holomorphic $k$-forms on $Z$} defined as the sheafification of the combinatorial sheaf $\SF^k$. With real coefficients, this can be defined as the kernel of the differential $d'' \colon \mathcal A^{k,0} \to \mathcal A^{k,1}$ from Dolbeault $(k,0)$-forms to Dolbeault $(k,1)$-forms on $Z$, as in~\cites{JSS19}. We give an alternative characterization of this sheaf that shows it is defined over $\Q$.

If $U$ is an open subset of $Z$, we say that $U$ is \emph{nice} if either $U$ is empty, or there exists a face $\gamma$ of $Z$ intersecting $U$ such that, for each face $\delta$ of $Z$, every connected component of $U\cap\delta$ contains $U\cap\gamma$. This implies that $U$ is connected, and that for each face $\delta\in Z$ which intersects $U$, $\gamma$ is a face of $\delta$ and the intersection $\delta\cap U$ is connected. We call $\gamma$ the \emph{minimum face of\/ $U$}. We have the following elementary result, whose proof is omitted.
\begin{prop}
  Nice open sets form a basis of open sets on $Z$.
\end{prop}

The sheaf $\ss\Omega^k_Z$ is the unique sheaf on $Z$ such that for each nice open set $U$ of $Z$ with minimum face $\gamma$, we have
\[ \ss\Omega^k_Z(U)=\SF^k(\gamma). \]

\medskip

Suppose now that $\comp Z$ is a compactification of $Z$ and denote by $\i\colon Z\hookrightarrow\comp Z$ the inclusion. We denote by $\ss\Omega^k_{Z,c}$ the \emph{sheaf of holomorphic $k$-forms on $\comp Z$ with compact support in $Z$} defined on connected open sets by
\[ \ss\Omega^k_{Z,c}(U):=\begin{cases}
\ss\Omega^k_Z(U) & \text{if $U\subseteq Z$,} \\
0             & \text{otherwise.}
\end{cases} \]

In other words, we have
\[ \ss\Omega^k_{Z,c}=\i_!\ss\Omega^k_Z, \]
that is, \emph{direct image with compact support} of the sheaf $\ss\Omega^k_Z$. The cohomology with compact support of $\ss\Omega^k_Z$ is computed by the usual sheaf cohomology of $\ss\Omega^k_{Z,c}$, \ie, we have
\[ H_c^\bul(Z, \ss\Omega^k_Z)=H^\bul(\comp Z, \ss\Omega^k_{Z,c}). \qedhere \]

\subsection{Cohomology with coefficients in the sheaf of tropical holomorphic forms}

Let $\Sigma$ be a simplicial homologically smooth tropical fan. For each cone $\sigma \in \Sigma$, we get the sheaf $\ss\Omega^k_{\compSigma^\sigma_\infty}$ of holomorphic $k$-forms on $\compSigma^\sigma_\infty \hookrightarrow \comp \Sigma$ which by pushforward leads to a sheaf on $\comp \Sigma$. We denote this sheaf by $\ss\Omega^k_\sigma$. The following proposition describes the cohomology of these sheaves.

\begin{prop} \label{prop:cohomology_Omega_fan} Notation as above, for each pair of non-negative integers $m, k$, we have
\[ H^m(\compSigma, \ss\Omega^k_\sigma) = \begin{cases}
  H^{2k}(\compSigma^\sigma) = H^{k,k}(\compSigma^\sigma) & \textrm{if $m = k$,} \\
  0 & \textrm{otherwise.}
\end{cases} \]
\end{prop}

\begin{proof}
We have
\[ H^m(\compSigma, \Omega^k_\sigma) \simeq H^m(\compSigma^\sigma_\infty, \Omega^k_{\compSigma^\sigma_\infty}) \simeq H^m(\compSigma^\sigma, \Omega^k_{\compSigma^\sigma}) \simeq H^{k,m}(\compSigma^\sigma). \]
The result now follows from Theorem~\ref{thm:ring_morphism_homology} and Proposition~\ref{prop:FY_duality}.
\end{proof}

\subsection{The resolution $\Omega^k_\bul$ of the sheaf $\Omega^k_{\Sigma,c}$} \label{sec:proof_resolution}

We will derive Theorem~\ref{thm:deligne-intro} by looking at the hypercohomology groups $\hyp^\bul(\comp\Sigma, \ss\Omega^k_\bul)$ of a resolution of $\Omega^k_{\Sigma,c}$.

For a pair of faces $\tau\subface\sigma$ in $\Sigma$, from the inclusion map $\compSigma^\sigma_\infty \hookrightarrow \compSigma^\tau_\infty \hookrightarrow \comp \Sigma$, we get a natural restriction map of sheaves $\i!_{\tau\subface\sigma}^*\colon\ss\Omega^k_\tau \to\ss\Omega^k_\sigma$ on $\comp \Sigma$. Here as before, the map $\i_{\tau\subface\sigma} = \i_{\sigma \supface\tau}$ denotes the inclusion $\compSigma^\sigma_\infty\hookrightarrow\compSigma^\tau_\infty$.

\smallskip
We consider now the following complex of sheaves on $\comp \Sigma$
\begin{equation}
\ss\Omega^k_\bul\colon \qquad \ss\Omega_\conezero^k \to \bigoplus_{\varrho \in \ssSigma_1}\ss\Omega^k_{\varrho} \to \bigoplus_{\sigma \in \ssSigma_2} \ss\Omega^k_{\sigma} \to \dots \to \bigoplus_{\sigma \in \ssSigma_{d-k}} \ss\Omega^k_{\sigma}
\end{equation}
concentrated in degrees $0, 1, \dots, d-k$, given by the dimension of the cones $\sigma$ in $\Sigma$. The boundary maps are given by
\[ \alpha\in\ss\Omega_\sigma^k \mapsto \d\alpha := \sum_{\zeta\ssupface\sigma}\sign(\sigma,\zeta)\i^*_{\sigma\ssubface\zeta}(\alpha). \]
Note that $\ss\Omega_{\Sigma, c}^k$ is a subsheaf of $\ss\Omega_\conezero^k$.

\begin{prop} \label{prop:exactness_Omega}
  The complex $\ss\Omega^k_\bul$ provides a resolution of $\ss\Omega_{\Sigma, c}^k$. That is, the following sequence of sheaves is exact
  \[ 0 \to \ss\Omega_{\Sigma, c}^k \to \ss\Omega_\conezero^k \to \bigoplus_{\varrho \in \ssSigma_1}\ss\Omega^k_\varrho \to \bigoplus_{\sigma \in \ssSigma_2}\ss\Omega^k_\sigma \to \dots \to \bigoplus_{\sigma \in \ssSigma_{d-k}}\ss\Omega^k_\sigma \to 0. \]
\end{prop}

\begin{proof}
  It will be enough to prove that taking sections over nice open sets $U$ give exact sequences of abelian groups.

  If $U$ is included in $\Sigma$, clearly by definition we have $\ss\Omega^k_{\Sigma,c}(U)\simeq\ss\Omega^k_{\conezero}(U)$, and the other sheaves of the sequence have no nontrivial section over $U$. Thus, the sequence is exact over $U$.

  It remains to prove that, for every nice open set $U$ having non-empty intersection with $\comp\Sigma\setminus\Sigma$, the sequence
  \[ 0 \to \ss\Omega_\conezero^k(U) \to \bigoplus_{\varrho \in \ssSigma_1}\ss\Omega^k_\varrho(U) \to \bigoplus_{\sigma \in \ssSigma_2}\ss\Omega^k_\sigma(U) \to \dots \to \bigoplus_{\sigma \in \ssSigma_{d-k}}\ss\Omega^k_\sigma(U) \to 0 \]
  is exact. Let $\gamma\in\comp\Sigma$ be the minimum face of $U$. Denote by $\sigma\in\Sigma$ the sedentarity of $\gamma$. The closed strata $\comp\Sigma_\infty^\tau \simeq \comp\Sigma^{\tau}$ of $\comp \Sigma$ which intersect $U$ are exactly those with $\tau\subface\sigma$. Moreover, if $\tau$ is a face of $\sigma$, we have
\[ \ss\Omega^k_\tau(U)=\SF^k(\gamma). \]
  Thus, the previous sequence can be rewritten in the form
  \begin{equation}\label{eq:simplex}
    0 \to \SF^k(\gamma) \to \bigoplus_{\tau\subface\sigma \\ \dims{\tau}=1}\SF^k(\gamma) \to \bigoplus_{\tau\subface\sigma \\ \dims{\tau}=2}\SF^k(\gamma) \to \dots \to \bigoplus_{\tau\subface\sigma \\ \dims{\tau}=\dims{\sigma}}\SF^k(\gamma) \to 0.
  \end{equation}
  This is just the cochain complex of the simplicial cohomology (for the natural simplicial structure induced by the faces) of the cone $\sigma$ with coefficients in the constant group $\SF^k(\gamma)$. This itself corresponds to the reduced simplicial cohomology of a simplex shifted by $1$. This last cohomology is trivial, thus the sequence is exact. That concludes the proof of the proposition.
\end{proof}

\begin{proof}[Proof of Theorem~\ref{thm:deligne-intro}]
  We have
  \[ H^m(\comp\Sigma,\ss\Omega^k_{\Sigma,c}) = H^m_c(\Sigma,\ss\Omega^k_\Sigma) = H^{k,m}_{c}(\Sigma)=
  \begin{cases}
    \SF^{d-k}(\conezero)^\dual = \SF_{d-k}(\conezero) & \text{if $m=d$,}\\
    0 & \text{otherwise.}
  \end{cases} \]
  Proposition~\ref{prop:exactness_Omega} implies that the cohomology of $\ss\Omega^k_{\Sigma,c}$ is isomorphic to the hypercohomology $\hyp(\Sigma, \ss\Omega^k_\bul)$. Therefore, we have
  \[ \hyp(\Sigma, \ss\Omega^k_\bul)\simeq \SF_{d-k}(\conezero)[-d], \]
  that is,
  \[ \hyp^m(\comp\Sigma, \ss\Omega^k_\bul) = \begin{cases}
    \SF_{d-k}(\conezero) & \textrm{for $m =d$,} \\
    0 & \textrm{otherwise.}
  \end{cases} \]

  On the other hand, the hypercohomology spectral sequence, combined with Proposition~\ref{prop:cohomology_Omega_fan}, implies that the hypercohomology of $\Omega^k_\bul$ is given by the cohomology of the following complex:
  \[ 0 \rightarrow H^{2k}(\comp\Sigma)[-k] \to \bigoplus_{\sigma \in \ssSigma_1} H^{2k}(\compSigma^\sigma)[-k-1] \to \dots \to \bigoplus_{\sigma \in \ssSigma_{d-k}} H^{2k}(\compSigma^\sigma)[-d] \to 0. \]

  Combining these two results, we conclude that the sequence
  \begin{equation}
  0 \rightarrow H^{2k}(\compSigma) \to \bigoplus_{\sigma \in \ssSigma_1} H^{2k}(\compSigma^\sigma) \to \bigoplus_{\sigma \in \ssSigma_2} H^{2k}(\compSigma^\sigma) \to \dots \to \bigoplus_{\sigma \in \ssSigma_{d-k}} H^{2k}(\compSigma^\sigma) \to \SF_{d-k}(\conezero) \to 0
  \end{equation}
  is exact, and the theorem follows.
\end{proof}

\subsection{An alternative proof using cellular cohomology} \label{sec:proof_cellular}

In this section we sketch another proof of Theorem \ref{thm:deligne-intro}. This proof uses the cellular cohomology we introduced in Section~\ref{sec:prel} instead of sheaf cohomology. The proof mainly follows the lines of the first proof. We simply states the spectral sequence analog to the sheaf spectral sequence, and the analogous of Proposition~\ref{prop:exactness_Omega}.

\begin{prop}
  The following sequence is exact for any integers $k$ and $q$.
  \[ 0 \to C_c^{k,q}(\Sigma) \to C^{k,q}(\comp\Sigma) \to \bigoplus_{\varrho \in \ssSigma_1}C^{k,q}(\compSigma^\varrho) \to \bigoplus_{\sigma \in \ssSigma_2}C^{k,q}(\compSigma^\sigma) \to \dots \to \bigoplus_{\sigma \in \ssSigma_{d-k}}C^{k,q}(\compSigma^\sigma) \to 0. \]
  The differentials are sums (with signs) of identity maps.
\end{prop}

\begin{proof}
  The sequence can be split into a direct sum for pairs $\tau\subface\eta$ of either
  \[ 0 \to \SF^k(\comp\eta) \to \SF^k(\comp\eta) \to 0 \]
  if $\tau=\conezero$, or
  \[ 0 \to 0 \to \SF^k(\cube_\eta^\tau) \to \bigoplus_{\varrho\in\Sigma_1\\\varrho\subface\tau} \SF^k(\cube_\eta^\tau) \to \bigoplus_{\sigma\in\Sigma_2\\\sigma\subface\tau} \SF^k(\cube_\eta^\tau) \to \dots \to \bigoplus_{\sigma\in\Sigma_{\dims\tau-1}\\\sigma\subface\tau} \SF^k(\cube_\eta^\tau) \to \SF^k(\cube_\eta^\tau) \to 0, \]
  otherwise. The summands for $\tau=\conezero$ are clearly exact. In the case $\tau\neq\conezero$, the sequence equals the tensor product of $\SF^k(\cube_\eta^\tau)$ with a sequence that can be identified we the standard Borel-Moore simplicial cochain complex of the cone $\tau$. The cohomology of this cochain complex is trivial, which concludes the proof.
\end{proof}

\begin{proof}[Sketch of an alternative proof of Theorem~\ref{thm:deligne-intro}]
  Theorem~\ref{thm:deligne-intro} now follows from a study similar to the one done in the previous proof, except that the spectral sequence is replaced by the one of the following double complex
  \[ \E^{a,b} := \begin{cases}
    C_c^{k,b}(\Sigma) & \text{if $a=-1$,} \\
    \bigoplus_{\sigma\in\ssSigma_a} C^{k,b}(\comp\Sigma^\sigma) & \text{otherwise.}
  \end{cases}. \]
  The horizontal differentials $\E^{\bul,\bul}\to \E^{\bul+1,\bul}$ are given by sums of identity maps. The vertical differentials $\E^{\bul,\bul}\to \E^{\bul,\bul+1}$ are the ones coming from tropical cohomology.
\end{proof}

\subsection{The case of general tropical fans} \label{subsec:Deligne_non_smooth}

In the case where $\Sigma$ is not necessarily homologically smooth, we still have the following proposition that we will use in Section~\ref{sec:homology_tropical_modification}.

\begin{prop} \label{prop:non_smooth_Deligne}
  Let $\Sigma$ be a simplicial tropical fan of dimension $d$, and let $k$ be an integer. The following sequence
  \[ \bigoplus_{\sigma \in \Sigma_{d-k-1}}H^{k,k}(\compSigma^\sigma) \longrightarrow \bigoplus_{\sigma \in \Sigma_{d-k}} H^{k,k}(\compSigma^\sigma) \longrightarrow H_c^{k,d}(\Sigma) \longrightarrow 0 \]
  is exact.
\end{prop}

Note that in the case $\Sigma$ is  homologically smooth, we have $H_c^{k,d}(\Sigma) \simeq \SF_{d-k}(\conezero)$, so in this case we recognize the end of the long exact sequence \eqref{eqn:Deligne_dual}.

\begin{proof}
  We proceed as in the previous case, and note that Proposition~\ref{prop:exactness_Omega} still holds in full generality, without requiring that $\Sigma$ be homologically smooth. This implies that the hypercohomology calculates the cohomology with compact support, that is,
  \[ \hyp^\bul(\Sigma,\Omega_\bul^k) \simeq H_c^{k,\bul}(\Sigma). \]

  Denote by $\E_\bul^{\bul,\bul}$ the hypercohomology spectral sequence. In page one, we have
  \[ \E_1^{a,b} = \bigoplus_{\sigma\in\Sigma_a} H^{k,b}(\comp\Sigma^\sigma). \]
  By Theorem~\ref{thm:ring_morphism_homology} applied to the star fan $\ssSigma^{\sigma}$, for $\sigma\in\Sigma$, we know that $H^{k,b}(\compSigma^\sigma)$ is trivial provided that $b > k$. Moreover, if $\dims\sigma > d-k$, then the dimension of $\ssSigma^\sigma$ is strictly less than $k$, which implies that $H^{k,b}(\compSigma^\sigma) = 0$. Therefore, $\E_1^{a,b}$ is nontrivial only for $0 \leq a \leq d-k$ and $0 \leq b \leq k$.

  Computing the further pages of the spectral sequence, we get that
  \[ \E_\infty^{d-k,k} = \E_2^{d-k,k} = \coker\bigl(\E_1^{d-k-1,k} \to \E_1^{d-k,k}\bigr). \]
  Since $\E_\infty^{d-k,k}$ is the only nontrivial term of degree $d$, we infer that
  \[ H_c^{k,d}(\Sigma) \simeq \E_\infty^{d-k,k} = \coker\Bigl(\bigoplus_{\sigma \in \Sigma_{d-k-1}}H^{k,k}(\compSigma^\sigma) \longrightarrow \bigoplus_{\sigma \in \Sigma_{d-k}} H^{k,k}(\compSigma^\sigma)\Bigr), \]
  which concludes the proof.
\end{proof}

\section{Coefficient sheaves of a tropical modification} \label{sec:coefficients_tropical_modification}

Let $(\Sigma,\omega!_\Sigma)$ be a homologically smooth tropical fan. Let $\Delta =\div(f)$ be a divisor associated to a meromorphic function $f$ on $\Sigma$, and let\/ $\~\Sigma = \tropmod{f}{\Sigma}$ be the tropical modification of\/ $\Sigma$ along $\Delta$. Let $\prtm \colon \~\Sigma \to \Sigma$ be the projection map.

In this section we prove the following theorem.
\begin{thm}[Local tropical modification formula] \label{thm:F_and_tropical_modification}
  Notation as above, for each non-negative integer $p$, the following holds.
  \begin{enumerate}[label=(\arabic*)]
    \item \label{thm:local_sigma}Let $\sigma$ be a face of\/ $\Sigma$ which is not included in $\Delta$. Then, the induced map $\prtm_*$ on coefficient sheaves gives an isomorphism
    \[
      \SF^{\~\Sigma}_p(\basetm\sigma) \simeq \SF^\Sigma_p(\sigma).
    \]
    \item \label{thm:local_delta}Let $\delta$ be a face in $ \Delta$.
    We have the following short exact sequences:
    \[ \begin{tikzcd}[row sep = tiny]
    0 \rar& \SF^\Delta_{p-1}(\delta) \rar& \SF^{\~\Sigma}_p(\uptm\delta) \rar{\prtm_*}& \SF^\Delta_p(\delta) \rar& 0, \\
    0 \rar& \SF^\Delta_{p-1}(\delta) \rar& \SF^{\~\Sigma}_p(\basetm\delta) \rar{\prtm_*}& \SF^\Sigma_p(\delta) \rar& 0.
    \end{tikzcd} \]
    The first map in both sequences is given by $\vv \mapsto \etm \wedge \prtm^*(\vv)$, with $\prtm^*(\vv)$ denoting an arbitrary preimage of\/ $\vv$ by the map $\prtm_*$ induced from the projection.
  \end{enumerate}
\end{thm}

Note that the map $\vv \mapsto \etm \wedge \prtm^*(\vv)$ does not depend on the choice of a preimage of $\vv$.

\begin{proof}
  We prove this theorem by using the tropical Deligne sequence.

  \smallskip
  For~\ref{thm:local_sigma}, let $\sigma$ be a cone in $\Sigma$ not included in $\Delta$.  Let $\ssf^\sigma$ be the meromorphic function induced by $f$ on $\ssSigma^\sigma$. Recall that this is obtained $\ssf^\sigma = \sspi^\sigma_*(f - \ell)$ for an element $\ell \in M$ that coincides with $f$ on $\sigma$, and $\sspi^\sigma\colon N \to \ssN^\sigma$ is the projection map.  Since $\sigma$ is not included in $\Delta = \div(f) = \div(f-\ell)$, there is no facet of $\Delta$ that contains $\sigma$. This means there is no codimension one face of $\Sigma$ that contains $\sigma$ and that is included in $\div(f)$. From this, we deduce that $\div(\ssf^\sigma)$ is trivial. Since $\Sigma^\sigma$ is homologically smooth, the Chow ring of $\ssSigma^\sigma$ verifies Poincaré duality. By characterization of div-faithfulness using the Chow rings given in~\cite[Section~4.5]{AP-hodge-fan}, we infer that $\ssSigma^\sigma$ is div-faithful.  As a consequence, given that $\div(\ssf^\sigma)=0$, we deduce that $\ssf^\sigma$ is induced by a linear map on $\ssN^\sigma$. By an abuse of the notation, we denote this linear map by $\ssf^\sigma$. Set now $\~\ell = \pi^{\sigma,*}(f^\sigma) + \ell$, and note that $\~\ell$ is a linear map that coincides with $f$ on the all faces of $\Sigma$ that contain $\sigma$. We infer that $\Gamma!_f$ and $\Gamma!_{\~\ell}$ coincide on each face of $\Sigma$ that contains $\sigma$. Moreover, $\Gamma!_{\~\ell}$ is a linear map which induces an isomorphism between $N_\R$ and $\Im(\Gamma!_{\~\ell})$. The isomorphism between $\SF^\Sigma_p(\sigma)$ and $\SF^{\~\Sigma}_p(\basetm\sigma)$ follows.

  \smallskip
  The first exact sequence in~\ref{thm:local_delta} is clear since $\uptm{\eta\/} \simeq \eta \times \ss\R_{\geq 0}\etm$ for any cone $\eta \supface \delta$ in $\Delta$.

  \smallskip
  It remains to prove the second assertion in~\ref{thm:local_delta}. We will prove the statement for $\delta =\conezero$, that is, the exactness of the following sequence
  \[ 0 \longrightarrow \SF^\Delta_{p-1}(\ss\conezero_\Delta) \longrightarrow \SF^{\~\Sigma}_p(\ss\conezero_{\~\Sigma}) \longrightarrow \SF^\Sigma_p(\ss\conezero_\Sigma) \longrightarrow 0. \]
  The statement for other faces $\delta \in \Delta$ then follows by applying the same argument to the tropical fans $\Delta^\delta$, $\ss{\~\Sigma}^{\basetm\delta}$, $\Sigma^\delta$, and using (non-canonical) compatible isomorphisms
  \[ \SF^\varUpsilon_k(\mu) \simeq \bigoplus_a \bigwedge^a N_\mu \,\otimes\, \SF^{\varUpsilon^\mu}_{k-a}(\ss\conezero_{\varUpsilon^\mu}) \]
  for each pair $(\mu,\varUpsilon)$ among $(\delta,\Delta), (\basetm\delta,\~\Sigma)$, or $(\delta,\Sigma)$ (the summands correspond to the graded pieces of the toric weight filtration on coefficient sheaves, see~\cite[Section 3.4]{AP-FY}).

  \smallskip
  Combining the sequence given by Proposition~\ref{prop:non_smooth_Deligne} for $\Delta, \Sigma$ and $\~\Sigma$, we get the diagram depicted in Figure~\ref{fig:dual_Deligne_tropical_modification} in which the three rows are exact.

  \begin{figure}[ht]
    \[ \begin{tikzcd}[cells={font=\everymath\expandafter{\the\everymath\displaystyle}}, row sep=scriptsize, column sep=scriptsize]
      0 \dar                                                                                 & 0 \dar \\
      \displaystyle\bigoplus_{\delta \in \Delta \\ \dims\delta=d-p-2}\!\!H^{p,p}(\compDelta^\delta)   \rar\dar& \bigoplus_{\delta \in \Delta \\ \dims\delta=d-p-1} H^{p,p}(\compDelta^\delta)     \rar\dar& H_c^{p,d-1}(\Delta) \rar\dar& 0 \\
      \bigoplus_{\eta \in \~\Sigma \\ \dims\eta=d-p-1}\!\!H^{p,p}(\ss{\comp{{\~\Sigma}}}^\eta) \rar\dar& \bigoplus_{\eta \in \~\Sigma \\ \dims\eta=d-p} H^{p,p}(\ss{\comp{{\~\Sigma}}}^\eta) \rar\dar& H_c^{p,d}(\~\Sigma) \rar\dar& 0 \\
      \bigoplus_{\sigma \in \Sigma \\ \dims\sigma=d-p-1}\!\!H^{p,p}(\compSigma^\sigma)   \rar\dar& \bigoplus_{\sigma \in \Sigma \\ \dims\sigma=d-p} H^{p,p}(\compSigma^\sigma)       \rar\dar& H_c^{p,d}(\Sigma)   \rar\dar& 0 \\
      0                                                                                      & 0                                                                                               & 0
    \end{tikzcd} \]
    \caption{Exact sequences given by Proposition~\ref{prop:non_smooth_Deligne} combined together for a tropical modification \label{fig:dual_Deligne_tropical_modification}}
  \end{figure}

  \smallskip
  We describe the vertical maps by providing the relation between cohomology groups of star fans and their compactifications, after tropical modification. Using Proposition~\ref{prop:star_fans_tropical_modifications}, for each face $\delta \in \Delta$, we get $\ss{{\~\Sigma}}^{\uptm\delta} \simeq \ss\Delta^\delta$. We thus obtain
  \[ H^{\bul,\bul}(\ss{{\~\Sigma}}^{\uptm\delta}) \simeq H^{\bul,\bul}(\ss\Delta^\delta) \qquad  \textrm{and} \qquad H^{\bul,\bul}(\ss{\overline{\~\Sigma}}^{\uptm\delta}) \simeq H^{\bul,\bul}(\ss{\compDelta}^\delta) \qquad \forall\, \delta \in \Delta. \]

  By the same proposition, if $\sigma$ is a face of $\Sigma$, we have $\ss{\~\Sigma}^{\basetm\sigma} \simeq \tropmod{\ssf^\sigma}{\Sigma^\sigma}$ where $\ssf^\sigma$ is the meromorphic function induced by $f$ on $\ssSigma^\sigma$. Since $\Sigma$ is homologically smooth, by Proposition~\ref{prop:FY_duality}, the Chow ring of $\ss\Sigma^\sigma$ verifies Poincaré duality. By the characterization of div-faithful tropical fans given in~\cite[Section 4.5]{AP-hodge-fan}, the fan $\ss\Sigma^\sigma$ is thus div-faithful. By the invariance of the Chow rings under tropical modifications for div-faithful tropical fans~\cite[Theorem 6.4]{AP-hodge-fan}, we deduce that
  \[ A^{\bul}(\ss{\~\Sigma}^{\basetm\sigma}) = A^{\bul}(\ss{\Sigma}^{\sigma}). \]
  Applying Theorem~\ref{thm:ring_morphism_homology} to the fans $\ss{\~\Sigma}^{\basetm\sigma}$ and $\ssSigma^{\sigma}$, we thus obtain
  \[ H^{p,p}(\ss{\comp{{\~\Sigma}}}^{{\basetm\sigma}}) \simeq H^{p,p}(\compSigma^\sigma). \]

  In all the cases, we have provided identifications of the cohomology of star fans and their compactifications in $\Sigma$ and $\~\Sigma$. Combining these isomorphisms for $\eta \in \~\Sigma$ of the form $\uptm\delta$ and $\basetm\sigma$, for $\delta \in\Delta$ and $\sigma\in\Sigma$, we get
  \begin{equation}\label{eq:split}
    \bigoplus_{\sigma\in\ss{\~\Sigma}_k} H^{p,p}(\ss{\comp{{\~\Sigma}}}^\sigma) = \bigoplus_{\delta \in \ss\Delta_{k-1}} H^{p,p}(\ss{\comp{{\~\Sigma}}}^{\uptm\delta}) \oplus \bigoplus_{\sigma \in \ss\Sigma_k} H^{p,p}(\ss{\comp{{\~\Sigma}}}^{\basetm\sigma}) \simeq \bigoplus_{\delta \in \ss\Delta_{k-1}} H^{p,p}(\ss{\comp{\Delta}}^{\delta}) \oplus \bigoplus_{\sigma \in \ss\Sigma_k} H^{p,p}(\ss{\comp{\Sigma}}^\sigma).
  \end{equation}
  The first two vertical short exact sequences in Figure~\ref{fig:dual_Deligne_tropical_modification} are split, and  obtained via isomorphisms given in~\eqref{eq:split} for $k=d-p-1$ and $k=d-p$. The maps in the last column are uniquely defined by the commutativity of the whole diagram. A diagram chasing now proves that the last column is also an exact sequence.

  \smallskip
  To conclude, we observe that we have a second commutative diagram, given as follows:
  \[ \begin{tikzcd}
        & H_c^{p,d-1}(\Delta)          \rar\dar[two heads]& H_c^{p,d}(\ss{\~\Sigma})                       \rar\dar[two heads]& H_c^{p,d}(\Sigma)        \rar\dar{\vsim}& 0 \\
  0 \rar& \SF^\Delta_{d-p-1}(\conezero_\Delta) \rar           & \SF^{\ss{\~\Sigma}}_{d-p}(\ss\conezero_{\~\Sigma}) \rar{\prtm_*}           & \SF^\Sigma_{d-p}(\ss\conezero_\Sigma) \rar& 0
  \end{tikzcd} \]
  The vertical maps are the ones given in~\eqref{eq:reformulation_PD_dual}, and they are all surjective, see Section~\ref{sec:tropical-fans}. Moreover, since $\Sigma$ is homologically smooth, it verifies tropical Poincaré duality, and so the last vertical map is an isomorphism. The injectivity of the map $\SF^\Delta_{d-p-1}(\ss\conezero_\Delta) \to \SF_{d-p}^{\~\Sigma}(\conezero_{\~\Sigma})$ is clear. We infer the exactness of the second row by a direct diagram chasing. This finishes the proof of Theorem~\ref{thm:F_and_tropical_modification}.
\end{proof}

\section{Homology of tropical modifications} \label{sec:homology_tropical_modification}

The aim of this section is to study the behavior of homology and cohomology groups under tropical modifications. Using this, we prove that homological smoothness for tropical fans is $\mT$-stable.

\subsection{Computation of the homology and cohomology of a tropical modification}

This section is devoted to the proof of the following theorem.

Let $(\Sigma,\omega!_\Sigma)$ be a homologically smooth tropical fan. Let $f$ be a meromorphic function on $\Sigma$ and let $(\Delta,\omega!_\Delta)$ be the divisor of $f$. Le $\~\Sigma = \tropmod{f}{\Sigma}$ be the tropical modification of $\Sigma$ along $\div(f)$. Let $\overline{\~\Sigma}$ be the canonical compactification of $\~\Sigma$. Let $\rho$ be the special ray $\uptm{\conezero\,}=\ss\R_{\geq 0} \uptm{\e}$ in $\~\Sigma$. Denote by $\ss{\~\Sigma}_{\infty}^{\rho}$ the corresponding stratum in the canonical compactification $\overline{\~\Sigma}$ of $\~\Sigma$. By Proposition~\ref{prop:star_fans_tropical_modifications}, this is isomorphic to $\Delta$. We consider $\~\Sigma\cup \ss{\~\Sigma}_{\infty}^{\rho}$ as an open subset of $\overline{\~\Sigma}$.

\begin{thm} \label{thm:homology_tropical_modification}
  Notation as above, the following equalities between different homology and cohomology groups hold.
  \begin{enumerate}[label=(\arabic*)]
    \item \label{thm:homology_tm_1}We have
    \[ H^{\bul,\bul}_c(\~\Sigma) \simeq H^{\bul,\bul}_c(\Sigma,\Delta) \qquad \textrm{and} \qquad
    H_{\bul,\bul}^\BM(\~\Sigma) \simeq H_{\bul,\bul}^\BM(\Sigma,\Delta). \]
    \item \label{thm:homology_tm_2}We have
    \[ H^{\bul,\bul}_c(\~\Sigma \cup \ss{\~\Sigma}^{\rho}_\infty) \simeq H^{\bul,\bul}_c(\Sigma)\qquad \textrm{and} \qquad H_{\bul,\bul}^\BM(\~\Sigma \cup \ss{\~\Sigma}^{\rho}_\infty)  \simeq H_{\bul,\bul}^\BM(\Sigma). \]

    \item \label{thm:homology_tm_3} We have
    \[ H^{\bul,\bul}(\overline{\~\Sigma}) \simeq H^{\bul,\bul}(\comp\Sigma) \qquad \textrm{and} \qquad H_{\bul,\bul}(\overline{\~\Sigma}) \simeq H_{\bul,\bul}(\comp\Sigma). \]
  \end{enumerate}
  These isomorphisms are all compatible with tropical Poincaré duality.
\end{thm}

We will describe  in the proof the different isomorphisms stated in the theorem. In a nutshell, these are induced  by the maps $\prtm_*\colon \SF_\bul(\basetm\sigma) \to \SF^\Sigma_\bul(\sigma)$ for $\sigma\in\Sigma$, and the inclusions $\i\colon \SF^{\ss{\~\Sigma}^{\rho}_\infty}_\bul(\delta) \simeq \SF^\Delta_\bul(\delta) \hookrightarrow \SF^\Sigma_\bul(\delta)$ for $\delta\in\Delta \simeq \ss{\~\Sigma}^{\rho}_\infty$ (this latter isomorphism follows from Proposition~\ref{prop:star_fans_tropical_modifications} and the description of strata in the canonical compactification~\ref{sec:canonical_compactification}).

Note that the notation $\SF_\bul(\delta)$ is ambiguous for $\delta\in\Delta$. So, we precise the fan in which  we consider $\SF_\bul(\delta)$ as a superscript: either $\SF^\Delta_\bul(\delta)$ or $\SF^\Sigma_\bul(\delta)$.

\begin{remark}
  In the above theorem, the open subset $\~\Sigma \cup \ss{\~\Sigma}^\rho_\infty$ can be viewed as the graph of $f$ in $\ssN_\R \times \T$, seen as a multivalued function over the tropical hyperfield. From the algebraic point of view, in the case $f$ is holomorphic, this is the analogue of the graph in $X\times \mathbb A^1$ of a holomorphic function $g$ defined on a subvariety $X$ of the algebraic torus. Clearly in the algebraic case, the cohomology does not change. The second statement in the theorem provides a tropical analogue of this assertion.
\end{remark}


As a first reduction, we note that the theorem is trivial in the case of a degenerate tropical modification, \ie, if $\Delta=0$. Indeed, since $\Sigma$ is homologically smooth, it is div-faithful by the characterization of div-faithful tropical fans given in~\cite[Section 4.5]{AP-hodge-fan} and so $f$ is linear. In this case,  $\tropmod{f}{\Sigma}$ is isomorphic to $\Sigma$, and the theorem follows easily. In what follows, we thus assume the tropical modification is non-degenerate.

We use the notation of Section~\ref{subsec:tropical_modification}, and denote by $\prtm$ the projection associated to the tropical modification. Also, $\etm$ is the unit vector of the special new ray $\rho$ in the tropical modification.

\smallskip
In the following, for two chain complexes $\left(C_\bul,\d\right)$ and $\left(D_\bul,\d\right)$ and a morphism $\phi\colon C_\bul \to D_\bul$, we denote by $\Cone_\bul(\phi)$ the \emph{mapping cone} chain complex defined by
\[ \Cone_k(\phi) \coloneqq C_{k-1} \oplus D_k,\quad\text{and}\quad \begin{array}{rccc}
  \partial\colon & C_{k-1} \oplus D_k & \longrightarrow & C_{k-2} \oplus D_{k-1},          \\
              &    a    \oplus  b  & \longmapsto     & -\d a   \oplus (\phi(a) + \d b).
\end{array} \]
We define an analogue notion of mapping cone cochain complex for a morphism of cochain complexes.  We refer to~\cite{KS90} for more details on mapping cones and their basic properties.

\begin{proof}
  We will use the description of the coefficients given in the previous section.

  We first prove the two isomorphisms stated in~\ref{thm:homology_tm_1}. Summing the short exact sequences of Theorem~\ref{thm:F_and_tropical_modification} over all faces of $\~\Sigma$, for each pair of integers $p$ and $q$, we get the following short exact sequence
  \[ \begin{tikzcd}[column sep=small]
    0  \rar&  C^\BM_{p-1,q-1}(\Delta) \oplus C^\BM_{p-1,q}(\Delta)  \rar&  C^\BM_{p,q}(\~\Sigma)  \rar&  C^\BM_{p,q-1}(\Delta) \oplus C^\BM_{p,q}(\Sigma)  \rar&  0.
  \end{tikzcd} \]

  An inspection of the boundary operators leads to the following short exact sequence
  \begin{equation} \label{eqn:short_exact_sequence_tropical_modification}
    \begin{tikzcd}[column sep=small]
      0  \rar&  \Cone_\bul\Bigl(C^\BM_{p-1,\bul}(\Delta) \xrightarrow{\id} C^\BM_{p-1,\bul}(\Delta)\Bigr)  \rar&  C^\BM_{p,\bul}(\~\Sigma)  \rar&  \Cone_\bul\Bigl(C^\BM_{p,\bul}(\Delta) \hookrightarrow C^\BM_{p,\bul}(\Sigma)\Bigr)  \rar&  0.
    \end{tikzcd}
  \end{equation}

  The homology of the cone of the identity map is always zero. Moreover, the homology of $\Cone_\bul\Bigl(C_{p,\bul}^\BM(\Delta) \hookrightarrow C_{p,\bul}^\BM(\Sigma)\Bigr)$ is isomorphic to the relative Borel-Moore homology of the pair $(\Sigma, \Delta)$, defined as the homology of the relative complex $C^\BM_{p,\bul}(\Sigma, \Delta)$.

  From the long exact sequence associated to the short exact sequence \eqref{eqn:short_exact_sequence_tropical_modification}, we get the following isomorphism
  \[ H^\BM_{\bul,\bul}(\~\Sigma)  \simeq  H^\BM_{\bul,\bul}(\Sigma, \Delta), \]
  which is exactly the first isomorphism of the theorem. The other isomorphism $H^{\bul,\bul}_c(\~\Sigma) \simeq H^{\bul,\bul}_c(\Sigma,\Delta)$ is obtained by duality.

  \medskip

  We now prove the isomorphism between $H_{\bul,\bul}^\BM(\~\Sigma \cup \ss{\~\Sigma}^{\rho}_\infty)$ and $H_{\bul,\bul}^\BM(\Sigma)$ stated in~\ref{thm:homology_tm_2}. Using the isomorphism between $\ss{\~\Sigma}^{\rho}_\infty$ and $\Delta$, the complex $C^\BM_{p,\bul}(\~\Sigma \cup \ss{\~\Sigma}^{\rho}_\infty)$ can be decomposed into the following short exact sequence
  \begin{equation} \label{eq:BM_relative_chain}
    \begin{split}
      0 \longrightarrow \Cone_\bul\Bigl(C^\BM_{p-1,\bul}(\Delta) \overset{\id_\Delta}{\longrightarrow} C^\BM_{p-1,\bul}(\Delta)\Bigr) \oplus C^\BM_{p,\bul}(\Delta) \longrightarrow C^\BM_{p,\bul}(\~\Sigma \cup \ss{\~\Sigma}^{\rho}_\infty) &\\
      & \hspace{-4cm} \longrightarrow \Cone_\bul\Bigl(C^\BM_{p,\bul}(\Delta) \overset{\i_\Delta}\longhookrightarrow C^\BM_{p,\bul}(\Sigma)\Bigr) \longrightarrow 0.
    \end{split}
  \end{equation}

  Comparing this short exact sequence with the short exact sequence of relative homology, we get the following commutative diagram:
  \[ \begin{tikzcd}[column sep = small, row sep = scriptsize]
    0 \rar& \Cone_\bul(\id_\Delta) \oplus C_{p,\bul}^\BM(\Delta) \dar{\pi_2}\rar& C_{p,\bul}^\BM(\~\Sigma \cup \ss{\~\Sigma}^{\rho}_\infty) \dar{\prtm_*+\i_{\Delta}}\rar& \Cone_\bul(\i_\Delta) \dar{\pi_1}\rar& 0\\
    0 \rar& C_{p,\bul}^\BM(\Delta) \rar& C_{p,\bul}^\BM(\Sigma) \rar& C_{p,\bul}^\BM(\Sigma,\Delta) \rar& 0
  \end{tikzcd} \]
  Here, $\ss\pi_2$ denotes the projection on the second part, $\prtm_*$ is the usual map coming from the projections $\prtm_*\colon \SF_p(\basetm\sigma) \to \SF^\Sigma_p(\sigma)$ for $\sigma\in\Sigma$. Moreover, using the isomorphism $\ss{\~\Sigma}^{\rho}_\infty \simeq \Delta$,  $\i!_{\Delta}$ is given by the maps $\SF_p^{\ss{\~\Sigma}^{\rho}_\infty}(\delta)\simeq\SF^\Delta_p(\delta)\hookrightarrow\SF^\Sigma_p(\delta)$ for $\delta\in\ss{\~\Sigma}^{\rho}_\infty \simeq \Delta$. The last map $\ss\pi_1$ is the natural projection
  \[ C_{p,\bul}^\BM(\Sigma) \oplus C_{p,\bul-1}^\BM(\Delta) \longrightarrow C_{p,\bul}^\BM(\Sigma) \longrightarrow \rquot{C_{p,\bul}^\BM(\Sigma)}{C_{p,\bul}^\BM(\Delta)}. \]

  We thus obtain a morphism between the long exact sequences associated to \eqref{eq:BM_relative_chain}:
  \[ \begin{tikzcd}[column sep = small, row sep = scriptsize]
    \cdots \rar& H^\BM_{p,q}(\Delta) \rar\dar[equal]& H^\BM_{p,q}(\~\Sigma \cup \ss{\~\Sigma}^{\rho}_\infty) \rar\dar{\prtm_*+\i_{\Delta}}& H_q(\Cone_\bul(\i_\Delta)) \rar\dar{\vsim}& H^\BM_{p,q-1}(\Delta) \rar\dar[equal]& \cdots \\
    \cdots \rar& H^\BM_{p,q}(\Delta) \rar& H^\BM_{p,q}(\Sigma) \rar& H^\BM_{p,q}(\Sigma,\Delta) \rar& H^\BM_{p,q-1}(\Delta) \rar& \cdots
  \end{tikzcd} \]
  The maps $\ss\pi_2$ and $\ss\pi_1$ induce isomorphisms in homology, and we conclude by using the five lemma. The other statement $H^{\bul,\bul}_c(\~\Sigma \cup \ss{\~\Sigma}^{\rho}_\infty) \simeq H^{\bul,\bul}_c(\Sigma)$ is proved by duality.

  \smallskip
  The last part of the theorem \ref{thm:homology_tm_3} can be proved in a similar way, we omit the details.
\end{proof}

\subsection{$\mT$-stability of homological smoothness}

As a consequence of Theorem~\ref{thm:homology_tropical_modification}, we prove Theorem~\ref{thm:smooth_shellable_intro}.

\begin{proof}[Proof of Theorem~\ref{thm:smooth_shellable_intro}]
  We use the $\mT$-stability meta lemma proved in~\cite{AP-hodge-fan}. It is a trivial fact that elements of $\Bsho$ are all homologically smooth. Moreover, since homological smoothness only depends on the support, we only have to prove the closedness by products and the closedness by tropical modifications. Closedness by products follows easily from Künneth formula. It remains to prove the closedness by tropical modifications.

  \smallskip
  Let $\Sigma$ be a tropical fan. Let $f$ be a conewise integral linear function on $\Sigma$. Set $\Delta = \div(f)$ and let $\~\Sigma = \tropmod{f}{\Sigma}$. Assume moreover that $\Delta$ is a non-empty homologically smooth tropical fan. By the meta lemma, we just need to prove tropical Poincaré duality for $\~\Sigma$.

  By Theorem~\ref{thm:homology_tropical_modification}, we know that $H_{\bul,\bul}^\BM(\~\Sigma) \simeq H^\BM_{\bul,\bul}(\Sigma,\Delta)$. Using the long exact sequence of relative homology, we get the following exact sequence
  \[ \dots \longrightarrow H^\BM_{p,q}(\Delta) \longrightarrow H^\BM_{p,q}(\Sigma) \longrightarrow H^\BM_{p,q}(\~\Sigma) \longrightarrow H^\BM_{p,q-1}(\Delta) \longrightarrow \cdots \]
  Since $\Sigma$ (\resp $\Delta$) verifies tropical Poincaré duality, its Borel-Moore homology is trivial except for $q=d$ (\resp $q=d-1$). We deduce that $H^\BM_{p,\bul}(\~\Sigma)$ is trivial except in degree $d$ and that we have a short exact sequence
  \[ 0 \longrightarrow H^\BM_{p,d}(\Sigma) \longrightarrow H^\BM_{p,d}(\~\Sigma) \longrightarrow H^\BM_{p,d-1}(\Delta) \longrightarrow 0. \]

  The cap products $\cdot \frown \nu_{\varUpsilon}$ for $\varUpsilon$ any of the tropical fans $\Sigma$, $\Delta$ or $\~\Sigma$ induces the following commutative diagram:
  \[ \begin{tikzcd}
  0 \rar& \SF^{d-p}(\conezero_\Sigma) \dar{\vsim}\rar& \SF^{d-p}(\conezero_{\~\Sigma}) \dar  \rar& \SF^{d-p-1}(\conezero_{\Delta}) \dar{\vsim}\rar& 0 \\
  0 \rar& H^\BM_{p,d}(\Sigma)     \rar& H^\BM_{p,d}(\~\Sigma)     \rar& H^\BM_{p,d-1}(\Delta)     \rar& 0
  \end{tikzcd} \]
  We have already seen that the second row is exact. The first row is also a short exact sequence by the dual of Theorem~\ref{thm:F_and_tropical_modification}. The first and last vertical maps are isomorphisms because $\Sigma$ and $\Delta$ both verify tropical Poincaré duality. Hence, by the five lemma, the second vertical map is an isomorphism, and therefore, $\~\Sigma$ verifies tropical Poincaré duality.
\end{proof}

\begin{remark} The theorem follows as well from Theorem~\ref{thm:homology_tropical_modification}, part~\ref{thm:homology_tm_3}, and the alternative characterization of tropical homology manifolds given in~\cite[Theorem~1.8]{AP-FY}.
\end{remark}

\section{Examples}\label{sec:examples}

In this final section, we consider some of the examples that appear in~\cite[Section 12]{AP-hodge-fan} and in~\cite[Section 8]{AP-FY} and expand them by focusing on properties that concern the tropical cohomology of tropical fans and their compactifications.

\subsection{The fan over the one-skeleton of the cube} \label{subsec:cube}

A rich source of examples is given by the fan defined over the one-skeleton of a cube. More precisely, consider the standard cube $\mbox{\mancube}$ with vertices $(\pm1, \pm1, \pm1)$, and let $\Sigma$ be the two-dimensional fan with rays generated by vertices and with facets generated by edges of the cube. The resulting fan $\Sigma$ endowed with orientation one on facets is tropical. Table~\ref{tab:cohomology_cube} summarizes the cohomological data of this fan.

\begin{table}
  \renewcommand{\strut}{\rule{0pt}{1.1em}}
  \[ \begin{array}{c|c|c|c|cc|c|c|c|}
  \cline{2-4} \cline{7-9}
  \strut H^{p,q}_c(\Sigma)       & q = 0 & 1    & 2                         & \qquad\qquad          & H^{p,q}(\comp\Sigma) & q = 0  & 1    & 2    \\ \cline{1-4} \cline{6-9}
  \multicolumn{1}{|c|}{\strut p = 0} & 0 & 0    & \Q^5                      & \multicolumn{1}{c|}{} & p = 0                    & \Q & 0    & 0    \\ \cline{1-4} \cline{6-9}
  \multicolumn{1}{|c|}{\strut 1} & 0 & 0    &                          \Q^3 & \multicolumn{1}{c|}{} & 1                    & 0  & \Q^5 & 0    \\ \cline{1-4} \cline{6-9}
  \multicolumn{1}{|c|}{\strut 2} & 0 & \Q^2 & \Q                        & \multicolumn{1}{c|}{} & 2                    & 0  & \Q^2 & \Q   \\ \cline{1-4} \cline{6-9}
  \end{array} \]
  \caption{Cohomology of the fan $\Sigma$ over the one-skeleton of the cube. \label{tab:cohomology_cube}}
\end{table}

Note that $\Sigma$ does not verify tropical Poincaré duality. However, using Theorem~\ref{thm:ring_morphism_homology}, it is possible to show that its Chow ring verifies Poincaré duality. The Chow ring of $\Sigma$ verifies the hard Lefschetz property and the Hodge-Riemann bilinear relations, so $\Sigma$ is Chow-Kähler, while it is not Kähler.

Moreover, the tropical Deligne sequence of $\Sigma$ for $p=2$
\[ 0 \to \SF^p(\conezero) \to \bigoplus_{\sigma \in \Sigma_2} H^0(\comp \Sigma^\sigma) \rightarrow \bigoplus_{\varrho \in \Sigma_1} H^2(\comp \Sigma^\varrho) \to H^{4} (\comp\Sigma) \to 0 \]
is not exact. Indeed, it can be easily checked that the Euler-Poincaré characteristic of the sequence is non-zero.

\subsection{The cross}

Consider the cross $\{xy = 0\}$ in $\R^2$ with four rays. We denote it by $\Delta$. It is the simplest tropical fan which is not homologically smooth.

\begin{example}[A tropical modification along a divisor which is not homologically smooth] \label{ex:trop_mod_along_non_smooth}
  Let $\Lambda$ be the complete fan in $\R^2$ with facets consisting of the four orthants (so that the one-skeleton is $\Delta$). It is easy to see that $\Delta$ is the divisor of the conewise integral linear function $f(x,y) = \min(0,x)+\min(0,y)$ on $\Lambda$. Consider $\Sigma\coloneqq \tropmod{f}{\Lambda}$. Since $\Lambda$ is homologically smooth, we get $H^{\bul,\bul}(\overline{\Sigma}) \simeq H^{\bul,\bul}(\comp{\Lambda})$. Hence, the cohomology of $\overline{\Sigma}$ verifies Poincaré duality. However, $\Sigma$ does not verify tropical Poincaré duality. Indeed, we have $H^\BM_{0,2}(\Sigma) \simeq \Q^4$ but $\SF^2(\conezero)\simeq \Q^3$.

  This shows that Theorem~\ref{thm:tropical_modification_intro} does not hold in general without the assumption of $\Delta$ being homologically smooth.
\end{example}

\begin{example}[A tropical modification which does not verify the local tropical modification formula]
  Consider the function $g$ on the cross $\Delta$ induced by $(x,y) \mapsto \max(0,x) - \max(0,y)$. Then, $\div(g)$ is trivial, and $\~\Delta \coloneqq \tropmod{g}{\Delta}$ is a tropical line in $\R^3$. Hence, $\SF_1^{\~\Delta}(\conezero) \simeq \Q^3$ is not isomorphic to $\SF_1^{\Delta}(\conezero) \simeq \Q^2$, whereas $\conezero\notin\div(g)$.

  This shows that Theorem~\ref{thm:F_and_tropical_modification} does not hold in general if $\Delta$ is not homologically smooth.
\end{example}

\begin{example}[A tropical modification along a nontrivial divisor which does not verify the local tropical modification formula]
  Consider the same function $g(x,y) = \max(0,x) - \max(0,y)$ as above on $\Lambda$. Let $\Sigma$ be the tropical modification introduced in Example~\ref{ex:trop_mod_along_non_smooth}, $\rho \coloneqq \uptm{\conezero\,} = \R_{\geq0}\etm \in \Sigma$ be the special new ray above $\conezero$ and $\prtm$ be the induced projection. Then, $\div(\prtm^*(g))$ is nontrivial: its support is the one-skeleton of $\Sigma$ deprived from $\rho$. Moreover, in the tropical modification $\~\Sigma \coloneqq \tropmod{\prtm^*(g)}{\Sigma}$, we have $\SF^{\~\Sigma}_1(\basetm\rho) \simeq \Q^4$ whereas $\SF^{\Sigma}_1(\rho) \simeq \Q^3$. So Theorem~\ref{thm:F_and_tropical_modification} does not hold.
\end{example}

\subsection{A quasilinear tropical fan which is not Bergman}

In~\cite[Example~12.17]{AP-hodge-fan}, we produce quasilinear  fans $\Sigma$ whose fanfolds are not Bergman. Theorem~\ref{thm:smooth_shellable_intro} implies that $\Sigma$ are all homologically smooth. Other examples of quasilinear fans appear in~\cite{Sch21}.

\subsection{The example from~\cite{BH17}}

An interesting example of a tropical fan was discovered by Babaee and Huh~\cite{BH17} in order to disprove a generalized version of the Hodge conjecture suggested in a work by Demailly. In~\cite[Example 12.24]{AP-hodge-fan} we study this examples in detail.  We mention here that this tropical fan is homologically smooth, but it is neither quasilinear nor Kähler.

\printbibliography

\end{document}